\documentclass{tMAM2e}% Use this class file. Do not change this line. This class file must be in your working directory.

\usepackage{epstopdf}% To incorporate .eps illustrations using PDFLaTeX, etc.
\usepackage[bookmarksnumbered,colorlinks,bookmarks,citecolor=blue,linkcolor=blue,urlcolor=blue,breaklinks,linktocpage]{hyperref}% To create hyperlinks using PDFLaTeX, etc.
\usepackage[all]{hypcap}% This turns on the hyperlinks for \label{} and \ref{}, and others.
\usepackage{subfigure}% Support for small, `sub' figures and tables

\usepackage{array}
\usepackage{amsthm}
\usepackage{amsmath}
\usepackage{amssymb}
\usepackage{fixltx2e}
\usepackage{graphicx}
\usepackage[all]{xy}
\usepackage{xcolor}
\usepackage{colortbl}
\usepackage{pdflscape}
\usepackage{multirow}
\usepackage{musicography}
\usepackage{musixtex}
\usepackage{hyperref}
\usepackage{caption}
\definecolor{green}{RGB}{0, 180, 130}
\definecolor{lightgray}{gray}{0.87}
\definecolor{orange}{RGB}{255, 87, 51}
\definecolor{blue-magenta}{RGB}{180, 138, 217}
\definecolor{purplepizzazz}{RGB}{254, 78, 218}

\theoremstyle{plain}
\newtheorem{theorem}{Theorem}[section]

\theoremstyle{definition}
\newtheorem{definition}[theorem]{Definition}

\theoremstyle{remark}

\newtheorem{example}[theorem]{Example}

\usepackage{nameref}

\begin{document}

% The following two lines are for production department only. Do not change.
%\jvol{00} \jnum{00} \jyear{2015}
%\doi{insert doi link}

% If this is an Editorial by guest editors, uncomment the next line. If this is a research article, do not change anything,
% since research article is the default and we do not mark articles as ``research articles'' anymore
%\articletype{EDITORIAL}

% Replace the title below by your title. Only the first letter and proper names are capitalized.
\title{A framework for topological music analysis (TMA)}

% Write out full first names, do not use initial. Middle initial or middle name(s) may be included if wished.
% Include the comma before the "and" if there are 3 authors.
% In the affiliations, use the format: Department, University, City, Country. No other information is included, such as street.
% In the affiliations, as already done below use a semicolon and \\ to separate successive affiliations, no semicolon on the last affiliation.
% Please use the exact formatting as in the example below, only replace your information. If there are two authors, their names should be separated by ``and''
% The submission requires an anonymous version and a non-anonymous version. To make the anonymous PDF, simply
% comment out author and affiliation information using %
\author{\name{Alberto Alcal\'a-Alvarez \textsuperscript{\ddag}$^{\ast}$\thanks{$^\ast$Corresponding author. Email: albertoalcala@ciencias.unam.mx} and Pablo Padilla-Longoria \textsuperscript{\ddag}}
\affil{\textsuperscript{\ddag}Institute for Applied Mathematics, National University of Mexico}
%leave blank for production:
%\received{}
}

\maketitle

% Replace the following abstract with your abstract in the ``abstract'' environment.
% Citations in the abstract must be written out in the full in Chicago style in parentheses, no hyperlinks to bibliography in the abstract.
% If you have an Online Supplement, mention that in the last line of the abstract, and say what its contents are. Capitalize both words: Online Supplement.
\begin{abstract}
In the present article we describe and discuss a framework for applying
different topological data analysis (TDA) techniques to a music fragment
given as a score in traditional Western notation. We first consider
different sets of points in Euclidean spaces of different dimensions
that correspond to musical events in the score, and obtain their persistent
homology features. Then we introduce two families of simplicial complexes
that can be associated with chord sequences, and leverage homology to compute their salient features. Finally, we show
the results of applying the described methods to the analysis and
stylistic comparison of fragments from three \emph{Brandenburg Concertos}
by J.S. Bach and two \emph{Graffiti} by Mexican composer Armando Luna.
\end{abstract}

% Please provide five to ten keywords taken from terms used in your manuscript in the ``keywords'' environment below.
% Separate keywords by semicolons followed by a space, no semicolon after the last word. No ``and'' preceding the final word.
% Italicize all foreign words, here and throughout article. German nouns are also capitalized. Delete the examples below.

% If the article has any substantial mathematics or computers at all, please include several American Mathematical Society (AMS) Mathematics Subject Classification (MCS) or the Association for Computing Machinery (ACM) Computing Classification Scheme (CCS) words. Please use the newest versions, currently 2010 and 2012 respectively.
% The AMS MCS codes can be found at: http://www.ams.org/msc/msc2010.html
% The ACM CCS words can be found at: http://www.acm.org/about/class/class/2012
% Separate by semicolons followed by a space, no semicolon after the last one. Delete the example numbers below.
% 00A65 Mathematics and music
\begin{classcode}\textit{2020 Mathematics Subject Classification}: 00A65; 55N31
%The example classification codes listed above translate as:
%00A65 Mathematics and music
%55N31 Persistent homology and applications, topological data analysis
%\\
%\textit{2012 Computing Classification Scheme}: applied computing
\end{classcode}

% Section and subsection titles only have the first letter capitalized, the rest is lowercase, except of course proper names

\section*{Introduction.}

In this work we present several ways of treating data extracted from
a digital music score, and discuss the results of applying some tools
and techniques from algebraic topology (mainly simplicial homology)
to music analysis. Our motivation is to incorporate new
scopes and computational tools to music analysis, hoping they will
contribute in establishing a theoretical and practical framework suitable
for analyzing music in a wide variety of languages and styles (ideally,
in \emph{any} language or style). 

Usually, music
is analyzed through the lens of a very specific framework, such as the traditional Western tonal theory,
the jazz modal harmonic setting, the dodecaphonic technique, the classical
Indian music tradition, etc. Of course all of these provide valuable
and useful analytical techniques and terms to deal with the elements
and processes which occur in their respective musical systems. Nevertheless,
when trying to describe, within a single framework, musical objects
and phenomena found in a diversity of repertoire or musical cultures,
musicians and musicologists often find less complete and consistent
methodological resources. Thus, for certain analytical and musicological purposes, it is desirable
to work in a more general framework within which to speak about 
many types and styles of music, in equal terms, though this
might necessarily lead to losing some of the fine details given by
more particular analytical scopes. 

In recent years, there has been
a wide range of works pointing in this direction, especially coming
from applied mathematics, introducing general theoretical and methodological frameworks that include the use ofconcepts and tools from  very different  mathematical areas (see, for example: \citet{AnMusStrucPerf,MazzolaTopos,TymoczkoGeometry,AndreattaMethodes,GrPattMatchPostTon,LluisGruposTeoriaMus,PareyonSelfSimilarity,EstradaTeoriad1,DynTopTools,IdEvMusStyIPadilla} ). In search of contributing to
this task, we focus on some mathematical models and tools that seem
pertinent to describe data codified (or codable) in musical scores
(computationally we will constrain ourselves to working with fragments
written in traditional Western music notation, though this does not
restrict us to Western music tradition or repertoire exclusively).
Of particular interest to us is the idea of applying techniques of
data analysis to musical information, and for the present article
we work with techniques from topological data analysis (TDA). The application of 
algebraic topology concepts and techniques to describe harmonic structure in
music has been treated in several works, for example:  \citet{TopSpMus,TopMusData,CompVisMusStrucSimpComp,TopStruCompMusAn,TopApprMusRel, DynTopTools,BigoAndreattaFiltration,HomPersTimeSeries,TopPersMachLearnMIR,  TopPersConvNNAudioSignals}.

Through the application of these techniques, we seek to deal directly with the events noted in a score, without assuming any given system of relations between pitches or pitch class sets. In contrast with the classical harmonic or Schenkerian  analysis, there is no assumption of a particular or pre-established hierarchical harmonic nor formal system. Also, we may state that our scope differs from other methodologies which include some “geometrization” of musical data,  such as the Tonnetz and its generalizations. Unlike these models, we do not assume any fixed structure of chords, and do not deal with voice-leading or chord-generation processes. We do not properly deal with \emph{chords} as defined in traditional music theory. Instead, we part from vertical events, that is, sets of pitches sounding simultaneously (as encoded in the score), which may actually incorporate two or more chords overlapping (as seen under a particular analytical framework). These vertical events are determined by the appearance, disappearance or prolongation of notes in the score.

On one hand, we study the persistent homology (for a presentation on the subject, the reader may refer to the \nameref{sec:Appendix} of this paper) of different sets of points
formed from pitch and time data in a score. We consider vertical events both
with and without their rhythm and onset in the score. Next,
we propose and explore two different ways of constructing sequences
of simplicial complexes from chord sequences. One of these constructions coincides
with similar scopes in the way of looking at a chord as a simplex
on vertices corresponding to pitches or pitch classes (see, for example,
\citet{CompVisMusStrucSimpComp}), modeling chords consisting
of $n$ pitch classes as simplices on $n$ vertices. Yet a novelty introduced in this
article is another way to construct a simplicial complex that encodes not only pitch but also interval information. Simplicial complexes associated with sequences of chords have been treated in previous works, such as \citet{DynTopTools,HomPersTimeSeries}. We plot as barcodes and persistence diagrams (see the \nameref{sec:Appendix}) the homological features of all such simplicial complexes associated to fragments of different scores,
and compare the results,
using the bottleneck distance (see the \nameref{sec:Appendix}). Finally,
we summarize our results in several dendrograms, and present our conclusions.

For the sake of space, we assume the reader is familiar with the basic
post-tonal theory concepts (pitch and interval classes; normal form
or order, and interval vector of a chord; see, for example \citet{IntroPostTonStraus}).
We present the basics of simplicial and persistent homology at the end of the paper, in the \nameref{sec:Appendix}. For a deeper treatment of homology (homology groups, Betti numbers, Euler characteristic) the reader can refer to, for example, the classic books of \citet{HatcherAlgTop,RotmanIntroAlgTop}. Persistent homology concepts and methods (Vietoris-Rips complexes, barcodes and persistence diagrams, bottleneck distance) may be consulted in the survey by \citet{PersHomSurvey}.

\section*{}

\section{Definitions and methods.}

We consider music scores written in traditional Western notation.
We define \textbf{vertical events} in a music score as tuples containing
information of synchronous sounds, usually pitches or pitch classes,
possibly together with some other features, such as its onset, duration (rhythm), dynamics (loudness) or timbre (instrumentation).
In this work we focus on events given as tuples of pitch classes,
with or without their duration and onset.

A vertical event given only by pitch information will be called a
\textbf{chord}, and a chord formed by $n$ different pitch classes
will be referred to as an $\boldsymbol{n}\text{\textendash}$\textbf{chord}.
A \textbf{music fragment} is a (usually assumed finite) sequence (ordered
set) of vertical events. Given a fragment $\mathcal{M}=\{\mathbf{e}_{0},\,...\,,\mathbf{e}_{N}\}$,
the \textbf{interval of events }$[\mathbf{e}_{i},\mathbf{e}_{j}]$
is the sequence $\{\mathbf{e}_{i},\mathbf{e}_{i+1},\,...\,,\mathbf{e}_{j}\}$. 

This segmentation of events in the score does not take into account any other information other than their sequential order. In the case of “continuous” (smooth) music passages or textures, one could identify starting points without the need of counting beats, and smooth fade-in/fade-out elements by a gradual discrete approximation. In any case, a digital measurable score is not absolutely necessary to apply the present model (though it actually could be produced). One may directly define the sequence of events $e_0,e_1, \, . . . $  

Working on sequences of events indexed by their order of occurrence allows us to easily focus merely in harmonic changes, and, when necessary, take into account time information (duration and onset, encoded as separate coordinates of an event). This implies that there may be chords appearing more than once in $\mathcal{M}$, when they occur at different
times (except in cases such as the example below, when a distinguished coordinate such as the onset is included). We will always take into account the index of a chord or event
in the sequence being considered. Thus, technically we should write
$\mathcal{M}=\{(e_{0},0),(e_{1},1),\,...\,,(e_{N},N)\}$. However, when there may be no confusion we will write only $e_{i}$ as to lighten notation. This notation also attempts to show the possible time-dependence of the consecutive occurrence of two or more chords, and helps us grasp harmonic progressions independently of the rhythmic values involved. Simplices corresponding to repeated chords will only appear once in the complex of cumulative events. Also, in the embedding of chords (without rhythm/onset) in $\mathbb R ^n$, repeated chords are mapped to the same point in $\mathbb R ^n.$

Given a music fragment $\mathcal{M}=\{\mathbf{e}_{0},\,...\,,\mathbf{e}_{N}\}$,
we define $\mathcal{A}(\mathcal{M})=\{a_{0},a_{1},\,...\,,a_{N}\}$
as the sequence (ordered set) whose $a_{i}$ term is the chord corresponding
to the $i\text{\textendash}$th vertical event of $\mathcal{M}$,
$\mathbf{e}_{i}$. Chords in the set $\mathcal{A}(\mathcal{M})$ may
be expressed in different ways: by their common name (C Major, d minor,
F\musSharp \ diminished, etc.), as a tuple of pitch classes (for
example, in ascending order, or following the normal form of the chord),
or as a tuple of intervals or interval classes (for example each chord
may be expressed as a sequence of intervals starting at the lowest
pitch (as in figured bass notation), or as an interval vector as defined
in classic post-tonal theory). 

For now we restrict ourselves to the setting of pitches within the
twelve tone equal temperament and more specifically, their corresponding
pitch classes represented as elements of the set of integers modulo
$12$, $\mathbb{Z}_{12}=\{\bar{0},\bar{1},\,...\,,\overline{11}\}$. This assumption allows us to show how these techniques can work in a fairly standard setting. Yet, the same methods can be applied in any other tempered, non-tempered or microtonal setting, making the proper straightforward adjustments to definitions dependant on the number of pitches or pitch classes considered (for example, the  mappings defined in \ref{sub: PersHomMaps}). We may think of our particular setting as a projection of the one considering all $1200$ \emph{cents} in one octave. It is also worth mentioning that the possibility of dealing with microvalues of pitch and duration implies that the methodology described may be useful to analyze \emph{continuous} music, including electronic music, as we could analyze discrete, measurable transcriptions of continuous musical textures.

We point out that the current  algorithm used to parse digital scores does not capture certain musical notations, such as grace notes, slurs (prolongation and phrasing), glissandi (unless explicitly written), tempo, metre, measure, and in general all text indications such as dynamics, playing mode, expressive marks, etc.

In order to conduct our harmonic analysis, we remove from the analyzed scores all staves containing
unpitched percussion instruments\footnote{In general, notes representing unpitched sounds could be included in the analysis, for instance
assigning a numerical ``pitch'' value sufficiently distant from the ones representing
actual pitches.}. Time durations and positions will be expressed (with decimal values) in quarter notes 
(e.g., one eight rhythm $=0.5$, etc.). We point out this does not imply restricting ourselves to notes with onsets on integer quarter beats. The quarter note value is only taken as time unit; we could well establish the use of miliseconds, for instance. It is pertinent to point out that the specific measure units chosen either for time or pitch will determine the point clouds and simplicial complexes associated to musical events. Nevertheless, homeomorphic data encodings (change of units) lead to similar shape features under TDA analyses, being homology a topological invariant.

Let us exemplify these definitions:
\begin{example}
\label{exa: ALP - G.H. a F.J.H. cc.1-2}Let us consider the first
two measures (see figure \ref{fig: ALP - G.H. a F.J.H. cc. 1-2})
of Mexican composer Armando Luna\footnote{Armando Luna Ponce (1964-2015). Mexican composer born in the city
of Chihuahua. He studied at the National Conservatory of Music of
Mexico (where he later taught composition and music analysis) and
the Carnegie Mellon University in the U.S.A.. He mainly produced chamber
and symphonic music for acoustic instruments, in what he came to name
a \emph{ludic-eclectic-neorampageous} style. Many of his pieces take
the \emph{suite} structure of brief movements as a model, incorporating
Renaissance and Baroque dances from the European traidition as well
as many other genres from both academic and folk origins. Also, a
considerable amount of his works are \emph{hommages} dedicated to
different composers of the Western academic music pantheon, whose
language and style are synthesized and reinterpreted.}'s \emph{Graffiti Hommage to Franz Joseph Haydn} (\emph{G.H. to F.J.H.}),
belonging to his \emph{Graffiti}, a series of miniatures written for
ensemble in 2006 as musical hommages dedicated to different composers.
Each \emph{Graffiti }is built upon the motif corresponding to the
musical translation of the name of the composer to whom it is dedicated,
according to the following letter equivalences: 

\begin{center}
A=la, B=si\musFlat, C=do, D=re, E=mi, F=fa, G=sol, H=si, S=mi\musFlat. \footnote{The equivalences of letters A,B,...,H correspond to the usual German letter system for the notes of the diatonic scale. The equivalence S=mi\musFlat \ also comes from the German 'Es', which stands for E\musFlat \ and is pronounced like the letter 's'. It was used for example by Shostakovich to musically encode his name as the motif re-mi-do-si (D-S-C-H). See, for example \href{https://en.wikipedia.org/wiki/DSCH_motif}{DSCH motif in Wikipedia}.}
\end{center}

These pieces incorporate some of the most characteristic elements
in the styles of composers Johann Sebastian Bach, Franz Joseph Haydn,
Bela Bart\'ok, Dave Brubeck, and seven others, as seen and condensed
by Luna through his own language. You may listen to the recording
of \emph{Graffiti} by the Present Music ensemble \href{https://www.youtube.com/watch?v=rKEccPGQedM&list=PLpdxhnDJxFP5yYvoFF99fMkQ9ESC37aha&index=11}{here}.\\

\begin{figure}[btp]
\begin{centering}
\includegraphics{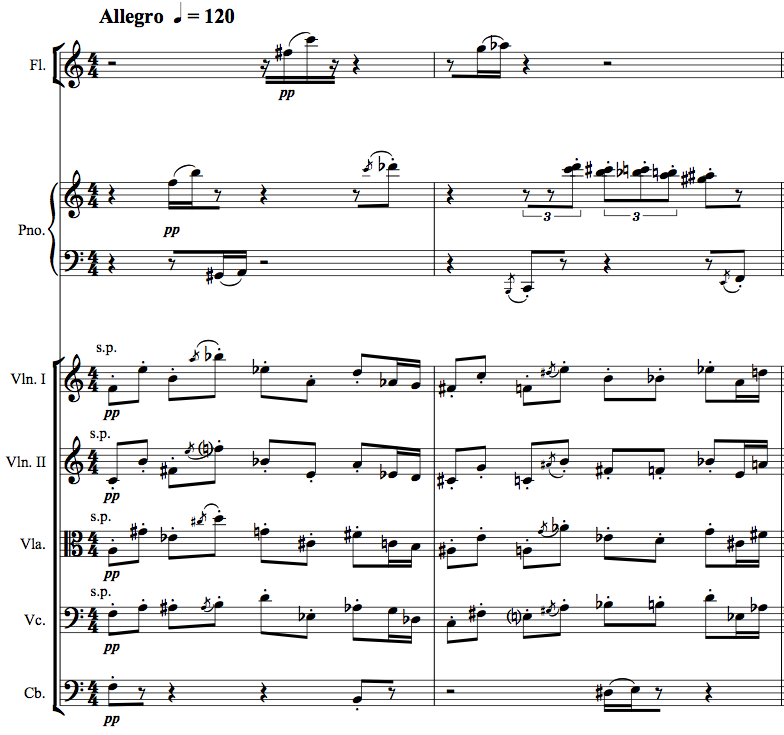}
\par\end{centering}

\protect\caption{{\label{fig: ALP - G.H. a F.J.H. cc. 1-2}Armando Luna-
\emph{G.H. to F.J.H.}, mm. 1-2. }}

\end{figure}

Below we present the sequence of the $28$ events of the music fragment
corresponding to the first two measures of Luna's \emph{G.H. to F.J.H},
given as triads of the form 
\[
(\text{normal form of chord},\text{ duration},\text{ onset})\,\text{.}
\]
That is, from score in figure \ref{fig: ALP - G.H. a F.J.H. cc. 1-2}
we get the following music fragment:

\begin{eqnarray*}
\mathcal{M} & = & \{\,((\bar{5},\bar{9},\bar{0}),0.5,0.0)\,,\,((\bar{4},\bar{8},\bar{9},\overline{11}),0.5,0.5),\\
 &  & ((\overline{10},\overline{11},\bar{3},\bar{5},\bar{6}),0.25,1.0)\,,\,((\overline{10},\overline{11},\bar{3},\bar{6}),0.25,1.25),\\
 &  & ((\bar{8},\bar{9},\overline{10},\overline{11},\bar{1},\bar{2},\bar{4},\bar{5}),0.25,1.5)\,,\,((\bar{9},\overline{10},\overline{11},\bar{2},\bar{5}),0.25,1.75),\\
 &  & ((\bar{2},\bar{3},\bar{7},\overline{10}),0.25,2.0)\,,\,((\bar{2},\bar{3},\bar{6},\bar{7},\overline{10}),0.25,2.25),\\
 &  & ((\bar{9},\bar{0},\bar{1},\bar{3},\bar{4}),0.25,2.5)\,,\,((\bar{9},\bar{1},\bar{3},\bar{4}),0.25,2.75),\\
 &  & ((\bar{6},\bar{8},\bar{9},\overline{11},\bar{2}),0.5,3.0)\,,\,((\bar{0},\bar{1},\bar{3},\bar{7},\bar{8}),0.25,3.5),\\
 &  & ((\bar{7},\overline{11},\bar{1},\bar{2}),0.25,3.75)\,,\,((\bar{6},\overline{10},\bar{0},\bar{1}),0.5,4.0),\\
 &  & ((\bar{0},\bar{4},\bar{6},\bar{7}),0.25,4.5)\,,\,((\bar{0},\bar{4},\bar{6},\bar{7},\bar{8}),0.25,4.75),\\
 &  & ((\bar{4},\bar{5},\bar{9},\overline{11},\bar{0}),0.333,5.0)\,,\,((\bar{4},\bar{5},\bar{9},\bar{0}),0.1666,5.333),\\
 &  & ((\bar{3},\bar{4},\bar{7},\bar{8},\bar{9},\overline{10},\overline{11}),0.1666,5.5)\,,\,((\bar{8},\bar{9},\overline{11},\bar{0},\bar{2},\bar{4}),0.333,5.666),\\
 &  & ((\overline{10},\overline{11},\bar{1},\bar{3},\bar{6}),0.25,6.0)\,,\,((\overline{10},\overline{11},\bar{1},\bar{3},\bar{4},\bar{6}),0.08333,6.25),\\
 &  & ((\overline{10},\overline{11},\bar{0},\bar{3},\bar{4},\bar{6}),0.1666,6.333)\,,\,((\overline{10},\overline{11},\bar{0},\bar{2},\bar{5}),0.1666,6.5),\\
 &  & ((\bar{9},\overline{10},\overline{11},\bar{2},\bar{5}),0.333,6.666)\,,\,((\bar{7},\bar{8},\overline{10},\bar{2},\bar{3}),0.5,7.0),\\
 &  & ((\bar{9},\bar{1},\bar{3},\bar{4},\bar{5}),0.25,7.5)\,,\,((\bar{2},\bar{5},\bar{6},\bar{8},\bar{9}),0.25,7.75)\,\}
\end{eqnarray*}
Of course we may get several music fragments for the same score, as
events may be described in several different ways, for example as
chords only, or adding some other relevant information found in the
digital score. Observe that an event is generated every time there
is a change in the notes sounding simultaneously according to the
score. Vertical events are taken as notes with the same rhythmic value.
Articulation (namely \emph{staccato}) does not affect the duration
of notes. If pitches last longer than the vertical event's associated
rhythmic value, some pitch classes will appear in the subsequent events
in order to fill their actual duration according to the score. Note
for example that events $9$ and $10$, respectively starting at time
positions $2.5$ and $2.75$, show the fact that chord $(\bar{9},\bar{0},\bar{1},\bar{3},\bar{4})$,
where pitch class $\bar{0}$ has rhythmic value $0.25$ (a sixteenth
note) in the score while the others last for $0.5$ quarters, is followed
by $(\bar{9},\bar{1},\bar{3},\bar{4})$, with a rhythmic value of
$0.25$. Besides that, grace notes are considered as part of the chord
they precede. We point out that the chord $(\bar{9},\overline{10},\overline{11},\bar{2},\bar{5})$
is the only repeated chord in this sequence, and it appears twice:
at the $6\text{th.}$ and then at the $25\text{th.}$ positions. Periodic
decimal temporal values have been truncated.

Now from the above fragment $\mathcal{M}$, we get
\begin{eqnarray*}
\mathcal{A}(\mathcal{M}) & = & \{\,(\bar{5},\bar{9},\bar{0})_0\,,\,(\bar{4},\bar{8},\bar{9},\overline{11})_1\,,\,(\overline{10},\overline{11},\bar{3},\bar{5},\bar{6})_2\,,\,(\overline{10},\overline{11},\bar{3},\bar{6})_3,\\
 &  & (\bar{8},\bar{9},\overline{10},\overline{11},\bar{1},\bar{2},\bar{4},\bar{5})_4\,,\,(\bar{9},\overline{10},\overline{11},\bar{2},\bar{5})_5\,,\,(\bar{2},\bar{3},\bar{7},\overline{10})_6\,,\,(\bar{2},\bar{3},\bar{6},\bar{7},\overline{10})_7,\\
 &  & (\bar{9},\bar{0},\bar{1},\bar{3},\bar{4})_8\,,\,(\bar{9},\bar{1},\bar{3},\bar{4})_9\,,\,(\bar{6},\bar{8},\bar{9},\overline{11},\bar{2})_{10}\,,\,(\bar{0},\bar{1},\bar{3},\bar{7},\bar{8})_{11}\,,\,(\bar{7},\overline{11},\bar{1},\bar{2})_{12},\\
 &  & (\bar{6},\overline{10},\bar{0},\bar{1})_{13}\,,\,(\bar{0},\bar{4},\bar{6},\bar{7})_{14}\,,\,(\bar{0},\bar{4},\bar{6},\bar{7},\bar{8})_{15}\,,\,(\bar{4},\bar{5},\bar{9},\overline{11},\bar{0})_{16}\,,\,(\bar{4},\bar{5},\bar{9},\bar{0})_{17},\\
 &  & (\bar{3},\bar{4},\bar{7},\bar{8},\bar{9},\overline{10},\overline{11})_{18}\,,\,(\bar{8},\bar{9},\overline{11},\bar{0},\bar{2},\bar{4})_{19}\,,\,(\overline{10},\overline{11},\bar{1},\bar{3},\bar{6})_{20}\,,\,(\overline{10},\overline{11},\bar{1},\bar{3},\bar{4},\bar{6})_{21},\\
 &  & (\overline{10},\overline{11},\bar{0},\bar{3},\bar{4},\bar{6})_{22}\,,\,(\overline{10},\overline{11},\bar{0},\bar{2},\bar{5})_{23}\,,\,(\bar{9},\overline{10},\overline{11},\bar{2},\bar{5})_{24}\,,\,(\bar{7},\bar{8},\overline{10},\bar{2},\bar{3})_{25},\\
 &  & (\bar{9},\bar{1},\bar{3},\bar{4},\bar{5})_{26}\,,\,(\bar{2},\bar{5},\bar{6},\bar{8},\bar{9})_{27}\,\}\,\text{.}
\end{eqnarray*}

\end{example}

Throughout the rest of the text, we will be ommiting subindices of chords and events.

In this work we deal with given events from different perspectives, encoding them in different ways, for example
\[
(\text{interval vector of chord},\text{ duration},\text{ onset})\,\text{.}
\]
Throughout this paper, we consider that events are given as tuples of the form
\[
(\text{chord as vector in some Euclidean space $\mathbb R ^n$},\text{duration},\text{onset})\,\text{,}
\]
or simply as vectors representing chords (without duration and onset). To deal with tuples constructed through a particular encoding of events, all those representing chords must belong to the same  $\mathbb R ^n$. In this article
we describe several ways of encoding chords as tuples in different
Euclidean spaces.

\subsection{General strategy.}

Given a music score, we extract its vertical events as tuples containing
the pitch classes in normal form, together with the duration
and onset in the score (both in quarter notes) of each
event. We also consider similar tuples in which chords are encoded as interval vectors (as defined in post-tonal theory). We point out that any of these representations of chords are a useful abstract standard for dealing with harmony in general (that is, outside the tonal context), but imply some loss of information, such as voicings and inversions of chords, as well as voice leading. To take into account such aspects, we need to consider some particular representation that encoded them, for instance modeling chord connections rather than chords themselves. This leads to other embeddings and associated spaces, such as the \textbf{connection simplicial complex} mentioned in the \nameref{sec:FutureWork} section. All such events will be
analyzed both with and without time values (duration and onset); that is to say, we will be analyzing chords with and without
their rhythm and position in time. 

One part of our analysis proposal involves classic TDA: musical events are encoded as
points in some Euclidean space $\mathbb{R}^{n}$ by means of different
embeddings (which we describe in the next subsection), to later calculate
their persistent homology under the Euclidean metric (using the Vietoris-Rips filtration; see \nameref{sec:Appendix}). 
On the other hand, we compute the homology of simplicial complexes directly associated with event intervals, without considering a metric among data points. We exhibit some cases in which these simplicial complexes actually form a filtration (though not associated with a metric).

All computations presented here were done in \href{https://www.python.org/}{Python}.
To parse and extract data from digital symbolic music files we use
the \href{http://web.mit.edu/music21/}{Music21} library. For constructing
simplicial complexes, obtaining their homologies and generating their
barcodes, persistence diagrams and the bottleneck distance between
them, we use algorithms from \href{https://mogutda.readthedocs.io}{MoguTDA},
Ripser and Persim (the latter two are incorporated as modules of the
\href{https://scikit-tda.org/}{Scikit-TDA} library). Also we made
use of some standard libraries for Mathematics and plotting, namely
NumPy, SciPy and Matplotlib. The scripts used for computations in this paper may be consulted in \href{https://github.com/al3-z/TMA}{this github project}.

We work with digital symbolic music files corresponding to musical
scores. The filetypes acceptable are those supported by Music21 and
include MIDI, .XML and .MXL files (Music21 supports many other \href{https://web.mit.edu/music21/doc/usersGuide/usersGuide_08_installingMusicXML.html}{formats}).
These digital scores are parsed in Python, after which different lists
of meaningful musical data are generated for a given fragment. Each
of these lists defines a musical data mapping. In our case, we obtain
lists of vertical events given as tuples consisting of the normal
form or interval vector of their corresponding chords, posssibly including
their duration (given in quarter notes) and offset (position in time
from the beginning of the score, also given in quarter notes). These
lists are the raw data analyzed by the persistent homology algorithms.
We generate and plot the corresponding persistence diagrams for different
mappings of data. As we said before, besides calculating the persistent homology of the Vietoris-Rips filtration constructed from a point cloud in some $\mathbb{R}^{n}$, we also compute and plot
the homological features of two families of simplicial complexes associated
with chord sequences in the score.

\subsection{Persistent homology on various musical data mappings. \label{sub: PersHomMaps}}

In this section we describe several sets of data points (in different
Euclidean spaces) associated with a given music fragment $\mathcal{M}$.
Concretely, we distinguish six different ways of generating tuples
in Euclidean spaces of different dimensions, from vertical events
in $\mathcal{M}$. We refer to these sets of tuples as data mappings
(as not to confuse them with \emph{data sets}, which usually would
mean the possible sets of scores or fragments analyzed), and number
them from I through VI. The first two mappings we consider (\colorbox{purplepizzazz}{\textbf{I}}
and \colorbox{cyan}{\textbf{II}}) consist of points that stand for
vertical events in a fragment $\mathcal{M}$ and incorporate temporal
data: each event is translated as a tuple of pitch or interval classes,
together with a rhythmic (duration) value and a onset
in time, both measured in quarter notes (see example \ref{exa: ALP - G.H. a F.J.H. cc.1-2}).
The rest of the mappings (\colorbox{yellow}{\textbf{III}}, \colorbox{orange}{\textbf{IV}},
\colorbox{green}{\textbf{V}}, and \colorbox{blue-magenta}{\textbf{VI}})
contain only pitch and/or interval information from each event, that
is, they focus only on harmonic aspects of $\mathcal{M}$. After generating
the data points of each mapping, we run upon each a persistent homology
analysis, and plot the resulting persitence diagrams or barcodes.
Then we compare persistence diagrams coming from different scores
in each mapping by using the bottleneck distance, which is a standard
tool for such task. Thus, we are able to depict in dendrograms these
distances among corresponding data mappings from such examples. This
lets us establish a certain  notion of \emph{closeness} between pieces
and styles. It is important to note that the persistent homology of
all mappings except I and V remains invariant under transposition
of the score by any interval; that is, the bottleneck distance between
persistence diagrams of a fragment and its transpositions is $0$
in all other mappings (II, III, IV and VI). So data mappings I and
V are the ones actually measuring the ``tonality'', ``tonic''
or current transposition of the score in question.

To illustrate the persistence diagrams obtained for each mapping,
we take the same fragment treated in example \ref{exa: ALP - G.H. a F.J.H. cc.1-2},
corresponding to the two measures shown in figure \ref{fig: ALP - G.H. a F.J.H. cc. 1-2}.
Colors in persistence diagrams correspond to different dimensions
(computed here up to dimension $3$, due to computational time \footnote{For some tests including higher dimensions in homology, ran over $32$-bar samples ($\sim 500$ event-points), there were overflow problems in the execution of the script.  As said in the \href{https://ripser.scikit-tda.org/en/latest/notebooks/Basic Usage.html}{documentation of the Ripser library}: ``[In] It [sic] practice, anything above $H_1$ is very slow.''}). For barcode plots we use two
different colors: teal for data mappings which include time information
(mappings I and II), and dark purple for the rest of mappings.
\begin{itemize}
\item \colorbox{purplepizzazz}{\textbf{Data mapping I}}: We begin considering
a set of points or vectors in $\mathbb{R}^{14}$ whose first twelve coordinates
correspond to pitch classes, followed by a rhythmic
and an onset values expressed in quarter notes. In these vectors, coordinates representing pitch classes
of a chord are represented by the integers $12,13,\,...\,,23$, according to their ordering in
the chord's normal form. We choose these representatives for each pitch class so we are able to embed
any $n$-chord, with $n=0,\, ...\,,11$ . The remaining
two coordinates (rhythm and onset) are given as a decimal
value representing length and onset in quarter notes.
For example, the first three events of the
score in figure \ref{fig: ALP - G.H. a F.J.H. cc. 1-2} yield the
following associated tuples in $\mathbb{R}^{14}$:
\begin{eqnarray*}
((\bar{5},\bar{9},\bar{0}),0.5,0) & \mapsto & (17,21,12,0,0,0,0,0,0,0,0,0,0.5,0)\\
((\bar{4},\bar{8},\bar{9},\overline{11}),0.5,0.5) & \mapsto & (16,20,21,23,0,0,0,0,0,0,0,0,0.5,0.5)\\
((\overline{10},\overline{11},\bar{3},\bar{5},\bar{6}),0.25,1) & \mapsto & (22,23,15,17,18,0,0,0,0,0,0,0,0.25,1)\,\text{.}
\end{eqnarray*}
Applying this mapping to the $28$-event fragment in figure \ref{fig: ALP - G.H. a F.J.H. cc. 1-2} (mm. 1-2 from Luna's \emph{Graffiti Hommage to F.J.H.}), we obtain the persistence and barcode diagrams shown
in figure \ref{fig:ALP-G.H. to F.J.H. cc1-2 - Persistence-and-barcode cat 1}.
\begin{figure}[btp]
\begin{centering}
\begin{minipage}[t]{0.45\columnwidth}%
\includegraphics[scale=0.5]{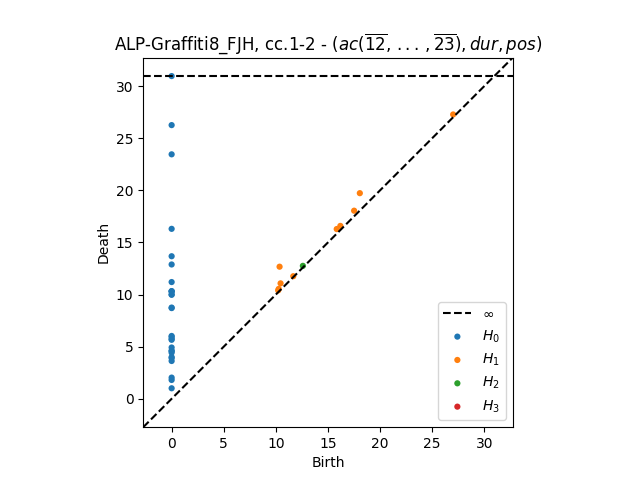}%
\end{minipage}\hfill{}%
\begin{minipage}[t]{0.45\columnwidth}%
\includegraphics[width=\textwidth,scale=0.5]{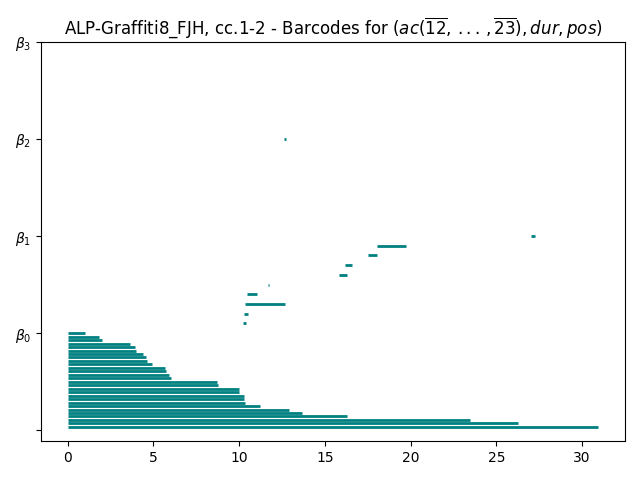}%
\end{minipage}
\par\end{centering}

\protect\caption{{\scriptsize{}\label{fig:ALP-G.H. to F.J.H. cc1-2 - Persistence-and-barcode cat 1}Persistence
and barcode diagrams from }\textbf{\scriptsize{}data mapping I}{\scriptsize{}
for mm. 1\textendash 2 of Luna's }\emph{\scriptsize{}Graffiti Hommage
to F.J.H.}.}
\end{figure}

Picturing events as pitch-rhythm-onset vectors, this mapping lets us get a notion of their general distribution over time, as well as identify the presence of distinguished harmonic regions. This way, together with an overview of the score, we can have a general impression of the harmonic-rhythmic texture of the fragment in question.

\item \colorbox{cyan}{\textbf{Data mapping II}}: Similarly, we consider
each event as a point in $\mathbb{R}^{8}$ whose first six coordinates
are the integers forming the interval vector of its corresponding
chord, followed by its duration and onset in time. This
mapping reflects similarity in the chord structures present in each
event, together with their distribution in time and rhythm. In this
case, for the same three events as above, we get:
\begin{eqnarray*}
((\bar{5},\bar{9},\bar{0}),0.5,0) & \mapsto & (0,0,1,1,1,0,0.5,0)\\
((\bar{4},\bar{8},\bar{9},\overline{11}),0.5,0.5) & \mapsto & (1,1,1,1,2,0,0.5,0.5)\\
((\overline{10},\overline{11},\bar{3},\bar{5},\bar{6}),0.25,1) & \mapsto & (2,1,1,2,3,1,0.25,1)\,\text{.}
\end{eqnarray*}
The resulting diagrams for example \ref{exa: ALP - G.H. a F.J.H. cc.1-2} are shown in figure
\ref{fig:ALP-G.H. to F.J.H. cc1-2 - Persistence-and-barcode cat 2}.
\begin{figure}[btp]
\begin{centering}
\begin{minipage}[t]{0.45\columnwidth}%
\includegraphics[scale=0.5]{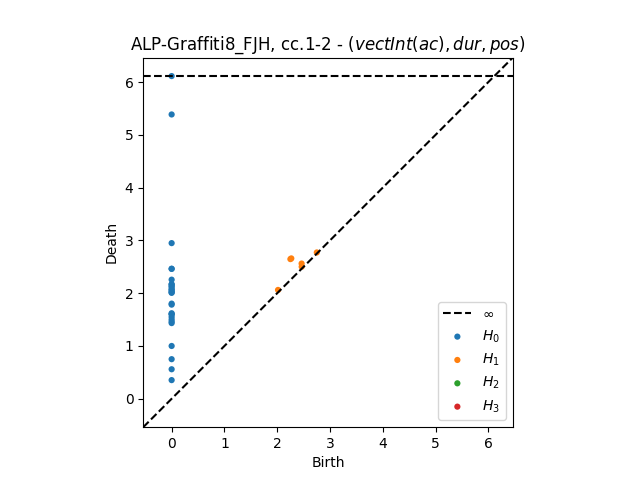}%
\end{minipage}\hfill{}%
\begin{minipage}[t]{0.45\columnwidth}%
\includegraphics[width=\textwidth,scale=0.5]{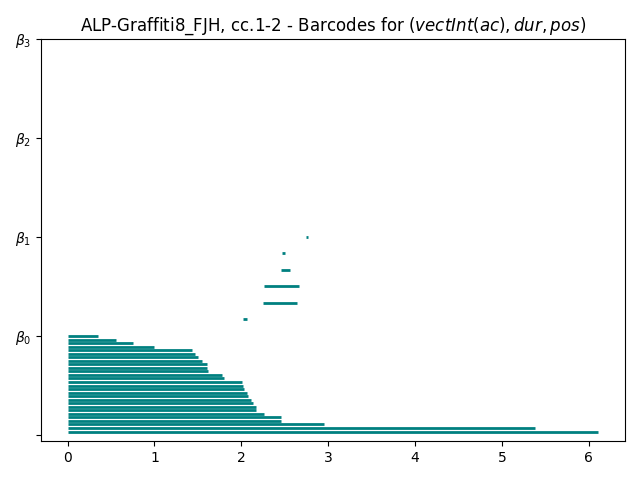}%
\end{minipage}
\par\end{centering}

\protect\caption{{\scriptsize{}\label{fig:ALP-G.H. to F.J.H. cc1-2 - Persistence-and-barcode cat 2}Persistence
and barcode diagrams from }\textbf{\scriptsize{}data mapping II}{\scriptsize{}
for mm. 1\textendash 2 of Luna's }\emph{\scriptsize{}Graffiti Hommage
to F.J.H.}.}
\end{figure}
This data mapping, together with mapping III focus on the types of
chords or different chord structures present in the fragment analyzed,
as they deal with interval vectors rather than the chords themselves.
\end{itemize}
Parallel to the above, we also work on sets obtained only from pitch
data. That is, we focus especially on harmony, by considering data
points containing only pitches in vertical events, forgetting about
their distribution along a timeline, their duration, etc. The points
generated for this analysis consist of tuples of pitch classes or
interval classes. We propose several ways of analyzing the same data,
by generating from a given fragment the following sets of points,
upon which we run a persistent homology analysis algorithm (under
the Euclidean distance):
\begin{itemize}
\item \colorbox{yellow}{\textbf{Data mapping III}}: Projection of data mapping II on its first six components. That is, we get vectors with integer
coordinates in $\mathbb{R}^{6}$, corresponding to interval vectors
of chords. Following example \ref{exa: ALP - G.H. a F.J.H. cc.1-2},
in this case we get the ordered set 
\begin{eqnarray*}
\mathcal{A}_{\text{int.vect.}}(\mathcal{M})= & \{\,(0,0,1,1,1,0)\,,\,(1,1,1,1,2,0)\,,\,(2,1,1,2,3,1)\,,\,(1,0,1,2,2,0),\\
 & (5,4,6,5,5,3)\,,\,(2,1,2,2,2,1)\,,\,(1,0,1,2,2,0)\,,\,(2,0,2,4,2,0),\\
 & (2,1,3,2,1,1)\,,\,(1,1,1,1,1,1)\,,\,(1,2,3,1,2,1)\,,\,(2,1,1,2,3,1),\\
 & (1,1,1,1,1,1)\,,\,(1,1,1,1,1,1)\,,\,(1,1,1,1,1,1)\,,\,(2,2,1,3,1,1),\\
 & (2,1,1,2,3,1)\,,\,(1,0,1,2,2,0)\,,\,(5,3,3,4,4,2)\,,\,(2,3,3,3,3,1),\\
 & (1,2,2,2,3,0)\,,\,(2,3,3,2,4,1)\,,\,(3,2,2,3,3,2),(2,2,2,1,2,1),\\
 & (2,1,2,2,2,1)\,,\,(2,1,1,2,3,1)\,,\,(2,2,1,3,1,1)\,,\,(2,1,3,2,1,1)\,\}\,\text{.}
\end{eqnarray*}

In figure \ref{fig:ALP-G.H. to F.J.H. cc1-2 - Persistence-and-barcode cat 3} we can see the persistence and barcode diagrams for the set $\mathcal{A}_{\text{int.vect.}}(\mathcal{M})$. In this case, they give us a hint of the sample's harmonic diversity.

This representation of events is the coarsest we consider. It is based purely on interval content modulo inversions (without duration or onset). This implies a simpler \emph{shape} of the point cloud, and so, less homological features present in its associated simplicial complex. Hence, persitence diagrams and barcodes for this mapping display less elements. Somehow they summarize diagrams obtained from all other mappings, as images of events under any of them can be ``projected'' onto the set of their interval vectors.

\begin{figure}[btp]
\begin{centering}
\begin{minipage}[t]{0.45\columnwidth}%
\includegraphics[scale=0.5]{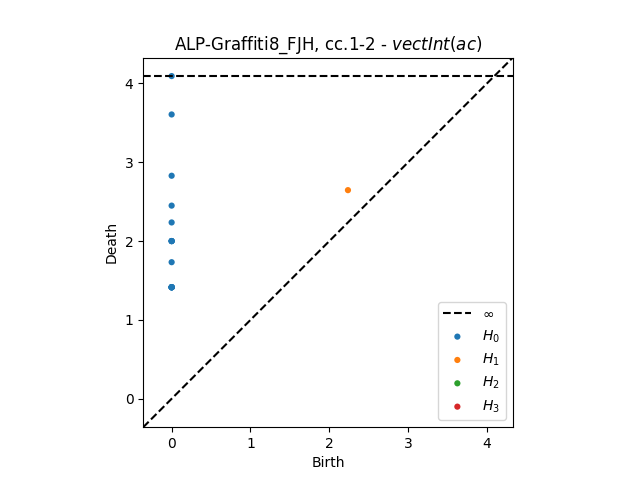}%
\end{minipage}\hfill{}%
\begin{minipage}[t]{0.45\columnwidth}%
\includegraphics[width=\textwidth,scale=0.5]{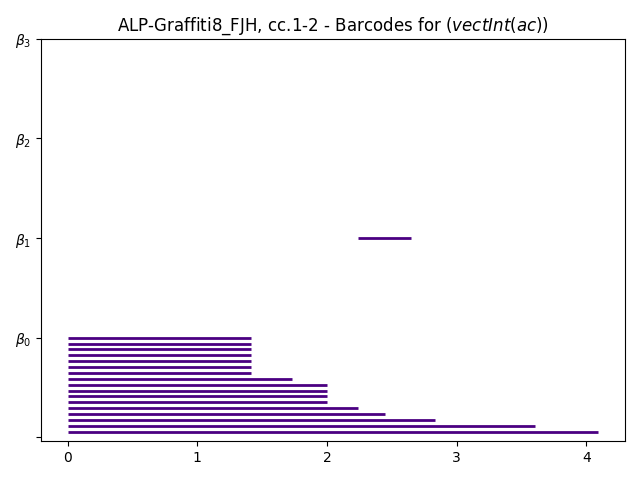}%
\end{minipage}
\par\end{centering}

\protect\caption{{\scriptsize{}\label{fig:ALP-G.H. to F.J.H. cc1-2 - Persistence-and-barcode cat 3}Persistence and barcode diagrams from }\textbf{\scriptsize{}data
mapping III}{\scriptsize{} for mm. 1\textendash 2 of Luna's }\emph{\scriptsize{}Graffiti
Hommage to F.J.H.}.}
\end{figure}

\item \colorbox{orange}{\textbf{Data mapping IV}}: Chords mapped as vectors
in $\{0,1\}^{12}\subset I^{12}\subset\mathbb{R}^{12}$, where $I=[0,1]\subset\mathbb{R}$,
as follows:

Given a chord $a$ in normal form, we define $a_{I^{12}}=(r_{0},\,...\,,r_{11})$
as:
\[
r_{i}=\begin{cases}
0 & \bar{i}\not\in a\\
1 & \bar{i}\in a
\end{cases}
\]
for $i\in\{0,1,\,...\,,11\}$ and $\bar{i}\in\{\bar{0},\bar{1},\,...\,,\overline{11}\}$,
the set of pitch classes represented by the smallest possible nonnegative
integers. For example, for the first chord of the score treated in
example \ref{exa: ALP - G.H. a F.J.H. cc.1-2}, $(\bar{5},\bar{9},\bar{0})$
(F Major vectorized following its normal form), we have
\[
(\bar{5},\bar{9},\bar{0})_{I^{12}}=(1,0,0,0,0,1,0,0,0,1,0,0,0)\,\text{.}
\]

This association yields the sequence 
\[
\mathcal{A}_{I^{12}}(\mathcal{M})=\{a_{i_{I^{12}}}\mid a_{i}\in\mathcal{A}_{\text{}}(\mathcal{M})\}\,\text{,}
\]
upon which we perform a persitent homology analysis, getting the
plots shown in figure \ref{fig:ALP-G.H. to F.J.H. cc1-2 - Persistence-and-barcode cat 4} for our particular example.

\begin{figure}[btp]
\begin{centering}
\begin{minipage}[t]{0.45\columnwidth}%
\includegraphics[scale=0.5]{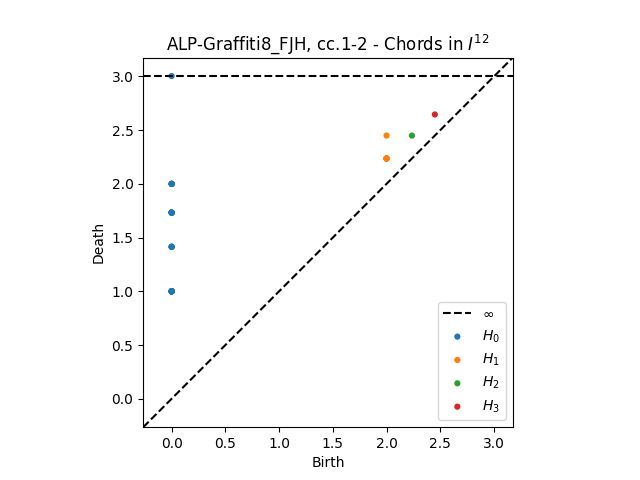}%
\end{minipage}\hfill{}%
\begin{minipage}[t]{0.45\columnwidth}%
\includegraphics[width=\textwidth,scale=0.5]{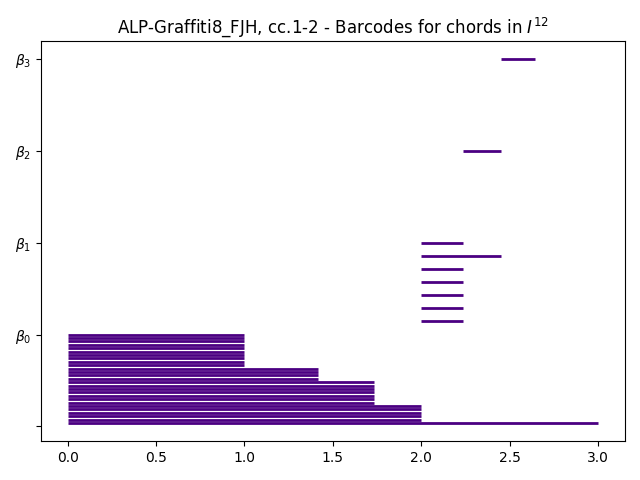}%
\end{minipage}
\par\end{centering}

\protect\caption{{\scriptsize{}\label{fig:ALP-G.H. to F.J.H. cc1-2 - Persistence-and-barcode cat 4}Persistence and barcode diagrams from }\textbf{\scriptsize{}data
mapping IV}{\scriptsize{} for mm. 1\textendash 2 of Luna's }\emph{\scriptsize{}Graffiti
Hommage to F.J.H.}.}
\end{figure}

In this setting, each dimension corresponds to a pitch class, and
so the images of two events are close exactly when they share most
of their pitches. More precisely, a given chord $a$ contains $n$
different pitch classes if and only if $\|a_{I^{12}}\|=\sqrt{n}$.
Furthermore, for chords $a,b$ we have $\|a_{I^{12}}-b_{I^{12}}\|=\sqrt{k}$
if and only if $a$ and $b$ differ in exactly $k$ pitch classes. This way, from the diagrams we can measure how close chords are among themselves, in terms of common/distinct pitches.

This mapping yields homological features in higher dimensions (up to $H_3$) than mappings III, V and VI. Also, these features persist only during specific intervals, determined by the square roots of integers. Thus, bars in barcodes appear to form ``blocks''.
More generally, we may take chords $a=\{x_{1},x_{2},\,...\,,x_{l}\}$
without octave equivalence of pitches. In this case, we may codify
each chord as a vector with integer coordinates in $\mathbb{R}^{12}$
taking this time
\[
r_{i}=\begin{cases}
0 & \bar{i}\not\in a\\
k & \bar{i}\,\,\text{appears exactly}\,\,k\,\,\text{times in}\,\,a
\end{cases}
\]

Note that mapping IV es the result of projecting the image of mapping VI onto the unit hypercube $I^{12}$, which explains the different shapes of their corresponding diagrams and barcodes.

\item \colorbox{green}{\textbf{Data mapping V}}: Projection of data mapping I on its first twelve components. In this case we get tuples of pitch classes as vectors in $\mathbb{R}^{12}$: for $n\in\{1,\,...\,,12\}$,
$n\text{\textendash}$chords are mapped to tuples with non-zero
integer values between $12$ and $23$ in the first $n$ entries, and $0$ in all the rest. As
a consequence, a chord is an $n\text{\textendash}$chord ($1\leq n\leq12$)
if and only if its associated vector belongs to the subspace spanned
by the first $n$ canonical basis vectors of $\mathbb{R}^{12}$, $\hat{\mathbf{e}}_{1},\,...\,,\hat{\mathbf{e}}_{n}$.
Thus, through this mapping, samples produce similar diagrams if and
only if their events are similar in pitch and number of harmonic voices.
For example, all triads belong to the $3\text{\textendash}$dimensional
linear space spanned by the first three standard basis vectors in
$\mathbb{R}^{12}$: 
\begin{eqnarray*}
\hat{\mathbf{e}}_{1} & = & (1,0,0,0\,...\,,0)\\
\hat{\mathbf{e}}_{2} & = & (0,1,0,0\,...\,,0)\\
\hat{\mathbf{e}}_{3} & = & (0,0,1,0\,...\,,0)\,\text{.}
\end{eqnarray*}
To associate such a tuple to a chord $a$, we choose integers $12,13,\,...\,,23$
as representatives of pitch classes $\bar{0},\bar{1},\,...\,,\overline{11}$,
and set each coordinate following the order given by the chord's normal
form, adding the necessary $0$s after the last pitch class. This
way, to chord $a=(\bar{5},\bar{9},\bar{0})$ we associate the tuple
\[
a(12,\,...\,,23)=a_{12^{+}}=(17,21,12,0,0,0,0,0,0,0,0,0)\in\mathbb{R}^{12}\,\text{.}
\]
We point out that in this mapping, chords are embedded into $\mathbb{R}^{12}$
as the vectors obtained from those in mapping I, without their last
two (time) coordinates.
\begin{figure}[btp]
\begin{centering}
\begin{minipage}[t]{0.45\columnwidth}%
\includegraphics[scale=0.5]{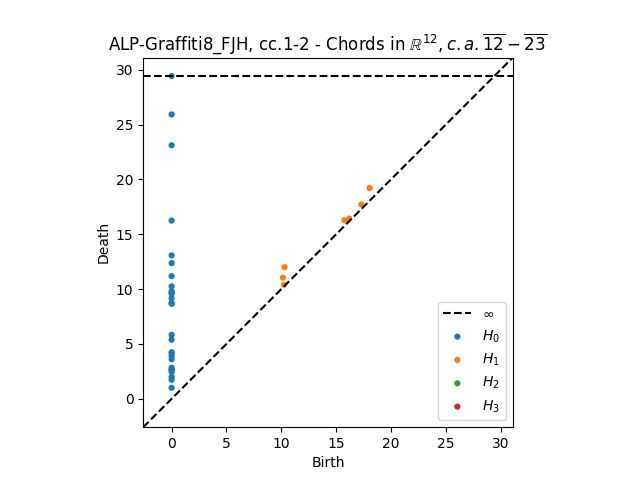}%
\end{minipage}\hfill{}%
\begin{minipage}[t]{0.45\columnwidth}%
\includegraphics[width=\textwidth,scale=0.5]{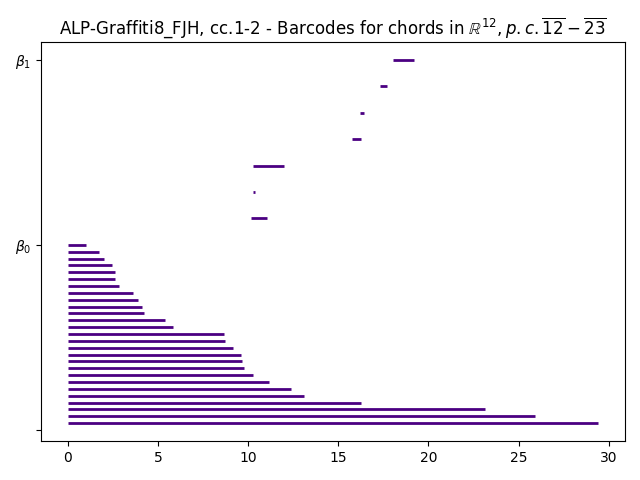}%
\end{minipage}
\par\end{centering}

\protect\caption{{\scriptsize{}\label{fig:ALP-G.H. to F.J.H. cc1-2 - Persistence-and-barcode cat 5}Persistence and barcode diagrams from }\textbf{\scriptsize{}data
mapping V}{\scriptsize{} for mm. 1\textendash 2 of Luna's }\emph{\scriptsize{}Graffiti
Hommage to F.J.H.}.}
\end{figure}
This mapping is one of the only two (the other one being mapping I, of which mapping V is a projection) that is sensitive
to transposition of the fragment in question by a given interval,
that is, the bottleneck distance between persistence diagrams for
this data mapping of a fragment and its transpositions is not always
$0$.

Through this mapping, we capture closeness of vertical events in terms of pitch content. As we said before, this mapping is sensitive to transpositions, due to the embedding of chords with $k$ pitch classes in the subspace generated by the first $k$ canonical basis vectors.
\item \colorbox{blue-magenta}{\textbf{Data mapping VI}}: Finally, we codify
not only the pitches in a chord, but also the intervals between two
consecutive pitch classes in a chord's normal form (or following some
other order for pitches). To do this, we propose the association $a\mapsto a_{\text{p},\text{int}}=(r_{0},\,...\,,r_{11})\in\mathbb{R}^{12}$,
where $0\leq r_{i}<12$ stands for the
interval class ($r_{i}>0$) between pitch class $\bar{i}$ and the
next pitch class in the normal form of $a$, if $\bar{i}$ belongs
to $a$, and $r_{i}=0$ otherwise. That is,
\[
r_{i}=\begin{cases}
0 & \bar{i}\not\in a\\
k & k\,\,\text{is the interval following}\,\,\bar{i}\,\,\text{in}\,\,a
\end{cases}
\]
In this case, we have, for chord $(\bar{5},\bar{9},\bar{0})$, vector
$(0,0,0,0,0,4,0,0,0,3,0,0)$. Under this mapping, points are close
to each other if and only if their corresponding events involve similar
intervals over the same pitches.
\begin{figure}[btp]
\begin{centering}
\begin{minipage}[t]{0.45\columnwidth}%
\includegraphics[scale=0.5]{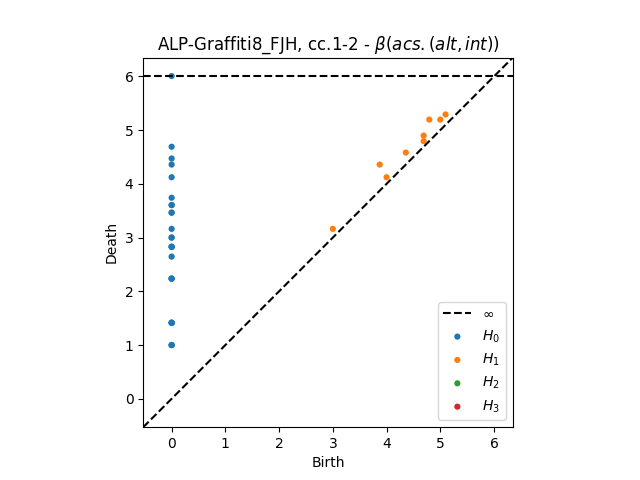}%
\end{minipage}\hfill{}%
\begin{minipage}[t]{0.45\columnwidth}%
\includegraphics[width=\textwidth,scale=0.5]{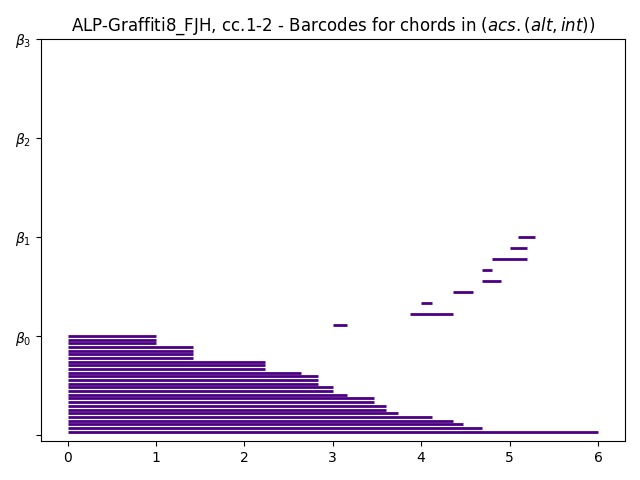}%
\end{minipage}
\par\end{centering}

\protect\caption{{\scriptsize{}\label{fig:ALP-G.H. to F.J.H. cc1-2 - Persistence-and-barcode cat 6}Persistence and barcode diagrams from }\textbf{\scriptsize{}data
mapping VI}{\scriptsize{} for mm. 1\textendash 2 of Luna's }\emph{\scriptsize{}Graffiti
Hommage to F.J.H.}.}
\end{figure}

Persistent homology analysis of data under this mapping usually shows non trivial homology cycles in higher dimensions than mappings III or V, sometimes agreeing with the dimensions of features detected by using mapping IV.

In contrast with mappings III and IV, bars in barcodes corresponding to mappings V and VI are more scattered. Thus, instead of forming ``blocks'', persistent features appear staggered.

\end{itemize}

Mappings IV, V and VI are different encodings of the normal form of chords. They are related to each other, since we may recover the pitch classes of a vertical event from any of them (and thus obtain their normal form). Nevertheless, they yield diagrams with different levels of detail and varying homological features. We hope that these variations let us have a more complete perspective on the harmonic data.

As we said before, we will use the resulting homological features
of our data mappings as stylistic descriptors. For instance, in the
case of the fragment we have chosen to illustrate our analysis proposal
(see example \ref{exa: ALP - G.H. a F.J.H. cc.1-2}), we identify
by sight a certain general pattern in the shapes of all the six pairs
of diagrams and barcodes: one or two connected components ($H_{0}$
/ $\beta_{0}$ values) that are present in most of Vietoris-Rips complexes
(see, for instance \citet{PersHomSurvey}) associated with the corresponding
cloud of data points, and a few briefly present hollowed circles ($H_{1}$
/ $\beta_{1}$ values). 

Following a standard TDA procedure, after calculating the persitent
homology of these six sets of points for different examples, we compute
the bottleneck distance (again, refer to \citet{PersHomSurvey}) between
their corresponding $H_{0}\text{\textendash}$diagrams, and plot a
dendrogram showing distances between them. We focus on $H_{0}\text{\textendash}$diagrams,
as only for a few samples and data mappings we obtained diagrams in
higher dimensions (nevertheless we include dendrograms for $H_{1}\text{\textendash}$diagrams
when available). For this discussion, see section \ref{sec:Results.}.

\subsection{\label{sub:Two-harmonic-simplicial-complexes}Two harmonic simplicial
complexes.}

We propose the construction of two different families of simplicial
complexes to describe the harmonic structure and evolution of a music
fragment. In a subsequent paper we will also develop the construction
of a familiy of simplicial complexes describing harmonic connections. These two constructions do not consider data to be embedded in some metric space, but result from assigning a simplicial complex to each vertical event in the score.
The simplicial complexes described here are formed of simplices
or simplicial complexes representing individual chords, which are
then combined into a bigger simplicial complex which somehow captures
the harmonic structure of the given interval of events. From all the
tests run so far\footnote{So far, we have run the present methods on over 100 fragments, from 13 classical, baroque, renaissance composers, as well as from traditional Indian and Mexican music.}, we can remark that for similar score samples, the associated complexes introduced here have similar Betti numbers. This will become more clear from the examples developed below. Thus,
these mathematical objects may be useful for musical style identification
and classification.

\subsubsection{\label{sub:Simplicial-complexes of chords by pitch}Simplicial complexes
of cumulative chords by pitch. }

This representation captures the pitches of chords as vertices of
simplices which are ``added'' together as events occur through the
score. 

Given an $(n+1)\text{\textendash}$chord $a$ with normal form vector
$(\overline{x_{0}},\overline{x_{1}},\,...\,,\overline{x_{n}})$, we
define its \textbf{associated} $\boldsymbol{n}\text{\textendash}$\textbf{simplex}
as $s(a)=\{\overline{x_{0}},\overline{x_{1}},\overline{x_{2}},\,...\,,\overline{x_{n}}\}$.
This association coincides with the one presented for example in \citet{TopSpMus,CompVisMusStrucSimpComp,TopStruCompMusAn}, though it is treated differently.
Given a music fragment $\mathcal{M}=\{e_{0},\,...\,,e_{N}\}$, we
consider its sequence of chords $\mathcal{A}(\mathcal{M})=\{a_{0},a_{1},\,...\,,a_{N}\}$,
which yields the sequence 
\[
s(a_{0}),s(a_{1}),\,...\,,s(a_{N})
\]
of associated simplices. For any integers $0\leq i\leq j\leq N$,
we define the \textbf{simplicial complex of cumulative chords by pitch}
in the interval of chords $[a_{i},a_{j}]$, denoted by $\mathcal{K}(i,j)$,
as the simplicial complex on simplices $s(a_{i}),s(a_{i+1}),\,...\,,s(a_{j})$,
together with all their faces. This is easily seen to comply with
the definition of a simplicial complex. Somehow, $\mathcal{K}(i,j)$
codifies the ``shape'' of the harmonic sequence or path from chord
$a_{i}$ to $a_{j}$. Nevertheless, this codification is not sensitive
to the order of appearance of chords. To take into account this order,
we propose the study of the following sequence of complexes $\mathcal{K}(i,j)$,
and some of its subsequences:
\begin{eqnarray*}
\mathcal{S}_{a}(\mathcal{M}) = \{\mathcal{K}(m,n)\}_{0\leq m,n \leq N} & = & \mathcal{K}(0,0),\mathcal{K}(0,1),\,...\,,\mathcal{K}(0,N),\mathcal{K}(1,1),\mathcal{K}(1,2),\,...\,,\mathcal{K}(1,N),\\
 &  & \mathcal{K}(2,2),\,...\,,\mathcal{K}(2,N),\,\,...\,,\mathcal{K}(N-1,N-1),\mathcal{K}(N-1,N), \\
 & & \mathcal{K}(N,N)\,\text{.}
\end{eqnarray*}

We have defined simplices and simplicial complexes from ordered sets of pitch classes as non-oriented objects. The order we have chosen on pitch class sets is given their normal form. However, this choice does not have an effect on the simplex by pitch associated to a chord: no matter the ordering, we get the abstract (non-oriented) simplex on the same set of pitch classes taken as vertices. However, the structure of the simplicial complex by pitches and intervals associated with a chord (defined below), will vary according to the order of pitches, thus yielding a different topological encoding of the same data. This could become useful when trying to preserve information of the actual intervals appearing in the score (chord voicing). In that case it may be congruent to drop the octave-equivalence hypothesis and work directly with pitches rather than pitch classes. On the other hand, we emphasize that each of the simplices and simplicial complexes associated with vertical events do not represent a sequence of notes, but a set of notes vertically coincident in the score. Oriented simplices and simplicial complexes associated to vertical events could be considered in this framework, being interpreted as encoding the position (voicing) of chords, from lowest to highest, for example.

For now, we work
with the main homological descriptors of each $\mathcal{K}(i,j)$:
its Betti numbers and Euler characteristic (refer to the \nameref{sec:Appendix}). Since in our context we
are considering chords consisting of up to $12$ equally tempered pitch classes, the maximum dimension
of the associated simplices and thus of the complexes $\mathcal{K}(i,j)$
is $11$, and so it suffices to compute their first twelve Betti numbers
$\beta_{0},\beta_{1},\,...\,,\beta_{11}$. We focus on the sequence
of simplicial complexes $\mathcal{K}(0,0),\mathcal{K}(0,1),\,...\,,\mathcal{K}(0,N)$,
which cover the full fragment, to get a picture of the change in the
topology of these accumulated successive harmonic events. Note that
with every step in this sequence, we add a simplicial complex to the
one we have so far built, namely in step $i$ we merge complex $\mathcal{K}(i,i)$
with $\mathcal{K}(0,i-1)$.

We now give an example to show how this simplicial complexes are built:
\begin{example}
Consider the fragment treated in example \ref{exa: ALP - G.H. a F.J.H. cc.1-2}.
Let us show how we build $\mathcal{K}(0,0)$, $\mathcal{K}(0,1)$
and $\mathcal{K}(0,2)$ from simplices on vertices in the set of pitch
classes $\{\bar{0},\,...\,,\overline{11}\}$. The first three chords
(according to their normal form) in  this fragment are
\[
(\bar{5},\bar{9},\bar{0})\,,\,(\bar{4},\bar{8},\bar{9},\overline{11})\,,\,(\overline{10},\overline{11},\bar{3},\bar{5},\bar{6})\,\text{.}
\]

$\mathcal{K}(0,0)$ is just the simplex $s((\bar{5},\bar{9},\bar{0}))=\{\bar{5},\bar{9},\bar{0}\}$,
which is the $2\text{\textendash}$simplex on the three vertices $\bar{5},\bar{9},\bar{0}$,
together with all its $1\text{\textendash}$ and $0\text{\textendash}$faces.
That is,

\[
\mathcal{K}(0,0)=\{\,\{\bar{0}\},\{\bar{5}\},\{\bar{9}\},\{\bar{0},\bar{5}\},\{\bar{0},\bar{9}\},\{\bar{5},\bar{9}\},\{\bar{5},\bar{9},\bar{0}\}\,\}\,\text{.}
\]

To visualize an abstract simplicial complex of any dimension (particularly
$>3$), we may draw a graph whose vertices are its $0\text{\textendash}$faces
(also called vertices) and whose edges are its $1\text{\textendash}$faces.
$2\text{\textendash}$faces will then be represented as closed $3\text{\textendash}$paths,
$3\text{\textendash}$faces as closed $4\text{\textendash}$paths, and so
on, though not every closed path in the graph will correspond
to a simplex (we use the term \emph{closed path} instead of \emph{cycle} to avoid confusion between homology cycles and cycles in a graph). As an example, for the simplex
described above we get the graph shown in figure \ref{fig: K(0,0)  ALP-G.H. a F.J.H. cc1-2}.
It is important to note that since this graph represents a $2\text{\textendash}$simplex
together with its faces, it must be interpreted as a full triangle,
i.e. vertices, perimeter, and area. We could get the same graph for
the sequence of chords $(\bar{0},\bar{5}),(\bar{0},\bar{9}),(\bar{5},\bar{9})$,
but in such a case we would have to picture it as the vertices and
perimeter only, without the triangle's inscribed area. So it is always
important to keep in mind what the picture of the graph is actually
representing. Another disadvantage of this graphic representation
is that sometimes when adding new simplices and their faces we may
get cycles which are actually voids in the simplicial complex, and
not faces (see figure \ref{fig: K(0,2)  ALP - G.H. a F.J.H. cc.1-2}).
Nevertheless, since simplices are grouped properly in the plot, we
still get a very good picture of our simplicial complexes.

\begin{figure}[btp]
\begin{centering}
\includegraphics[scale=0.5]{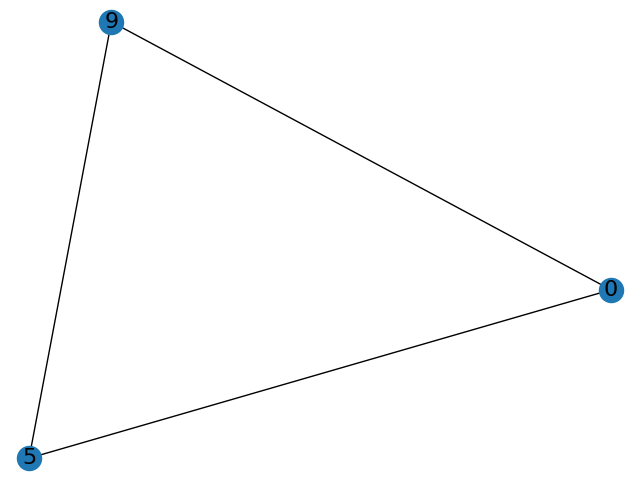}
\par\end{centering}

\protect\caption{{\label{fig: K(0,0)  ALP-G.H. a F.J.H. cc1-2}Simplicial}
complex $\mathcal{K}(0,0)$, corresponding to the chord $(\bar{5},\bar{9},\bar{0})$
(see example \ref{exa: ALP - G.H. a F.J.H. cc.1-2}).}

\end{figure}

$\mathcal{K}(0,1)$ is the simplicial complex whose simplices are
the simplices of $\mathcal{K}(0,0)$ together with simplex $s((\bar{4},\bar{8},\bar{9},\overline{11}))$
and its fifteen faces. So we get the following simplicial complex
on the six vertices $\bar{5},\bar{9},\bar{0},\bar{4},\bar{8},\overline{11}$:

\begin{eqnarray*}
\mathcal{K}(0,1) & = & \{\,\{\bar{0}\},\{\bar{5}\},\{\bar{9}\},\{\bar{4}\},\{\bar{8}\},\{\overline{11}\},\{\bar{0},\bar{5}\},\{\bar{0},\bar{9}\},\{\bar{5},\bar{9}\},\{\bar{4},\bar{8}\},\{\bar{4},\bar{9}\},\\
 &  & \,\,\{\bar{4},\overline{11}\},\{\bar{8},\bar{9}\},\{\bar{8},\overline{11}\},\{\bar{9},\overline{11}\},\{\bar{5},\bar{9},\bar{0}\},\{\bar{4},\bar{8},\bar{9}\},\\
 &  & \,\,\{\bar{4},\bar{8},\overline{11}\},\{\bar{4},\bar{9},\overline{11}\},\{\bar{8},\bar{9},\overline{11}\},\{\bar{4},\bar{8},\bar{9},\overline{11}\}\,\}\,\text{.}
\end{eqnarray*}
This complex has $21$ simplices: the seven simplices from $\mathcal{K}(0,0)$,
together with the fifteen faces of $s((\bar{4},\bar{8},\bar{9},\overline{11}))$,
out of which one is already in $\mathcal{K}(0,0)$ (the $0\text{\textendash}$face
$\{\bar{9}\}$). We get a picture of simplicial complex $\mathcal{K}(0,1)$,
corresponding to the sequence of chords $(\bar{5},\bar{9},\bar{0})\,,\,(\bar{4},\bar{8},\bar{9},\overline{11})$
in figure \ref{fig: K(0,1)  ALP-G.H. a F.J.H. cc1-2}. In this case,
we can clearly see the simplicial complex $\mathcal{K}(0,0)$, which
stands for chord $(\bar{5},\bar{9},\bar{0})$ (see figure \ref{fig: K(0,0)  ALP-G.H. a F.J.H. cc1-2}),
and the new added simplices forming simplicial complex $\mathcal{K}(1,1)$,
associated with chord $(\bar{4},\bar{8},\bar{9},\overline{11})$, which
again we must picture as a full tetrahedron containing all its vertices,
edges, faces, and volume. Both of these simplicial complexes are joined
together by pitch class $\bar{9}$, which is the only one common to
both chords.

\begin{figure}[btp]
\begin{centering}
\includegraphics[scale=0.5]{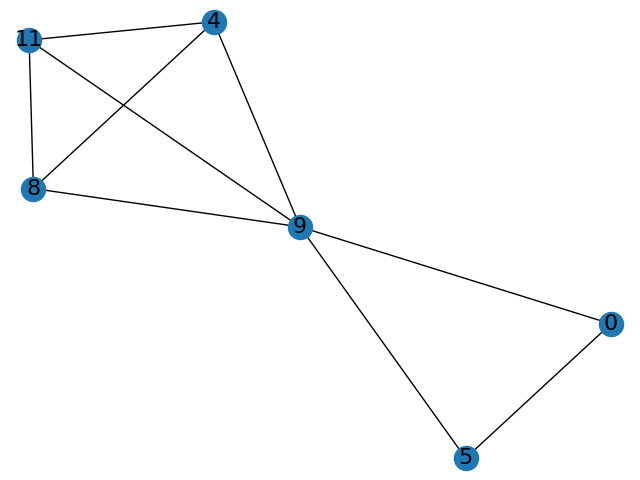}
\par\end{centering}

\protect\caption{{\label{fig: K(0,1)  ALP-G.H. a F.J.H. cc1-2}Simplicial}
complex $\mathcal{K}(0,1)$, corresponding to the sequence of chords
$(\bar{5},\bar{9},\bar{0})\,,\,(\bar{4},\bar{8},\bar{9},\overline{11})$
(see example \ref{exa: ALP - G.H. a F.J.H. cc.1-2}).}
\end{figure}

$\mathcal{K}(0,2)$ is in this case the simplicial complex formed
by all simplices in $\mathcal{K}(0,1)$, together with all $31$ faces
of simplex $\{\overline{10},\overline{11},\bar{3},\bar{5},\bar{6}\}$
(that accounts for five $0\text{\textendash}$faces, ten $1\text{\textendash}$faces,
ten $2\text{\textendash}$faces, five $3\text{\textendash}$faces
and one $4\text{\textendash}$face). From these, two have already
appeared in $\mathcal{K}(0,1)$, namely the two $0\text{\textendash}$faces
$\{\bar{5}\}$ and $\{\overline{11}\}$. So actually $\mathcal{K}(0,2)$
consists of $21+29=50$ simplices: nine $0\text{\textendash}$simplices,
nineteen $1\text{\textendash}$simplices, fifteen $2\text{\textendash}$simplices,
six $3\text{\textendash}$simplices, and one $4\text{\textendash}$simplex.

\begin{figure}[btp]
\begin{centering}
\includegraphics[scale=0.5]{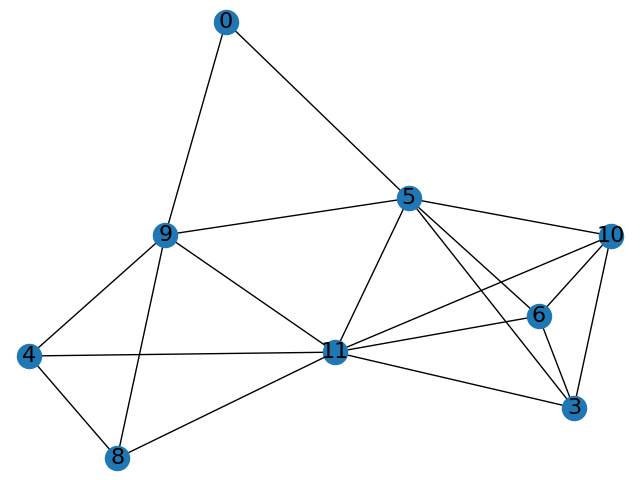}
\par\end{centering}

\protect\caption{\label{fig: K(0,2)  ALP - G.H. a F.J.H. cc.1-2}{Simplicial}
complex $\mathcal{K}(0,2)$, corresponding to the sequence of chords
$(\bar{5},\bar{9},\bar{0})\,,\,(\bar{4},\bar{8},\bar{9},\overline{11})\,,\,(\overline{10},\overline{11},\bar{3},\bar{5},\bar{6})$
(see example \ref{exa: ALP - G.H. a F.J.H. cc.1-2}).}
\end{figure}

We point out that in figure \ref{fig: K(0,2)  ALP - G.H. a F.J.H. cc.1-2}
we get a $3\text{\textendash}$cycle, the one formed by vertices labeled
$5,9,11$, which does not correspond to a simplex in $\mathcal{K}(0,2)$,
as it is not a face a simplex associated with any of the chords considered.
So actually this simplicial complex has a circular ($1\text{\textendash}$dimensional)
void given by this $3\text{\textendash}$cycle.

In figure \ref{fig:Barcodes K(0,i)  ALP-G.H. a F.J.H. cc1-2} we show
the barcode plot of Betti numbers for simplicial complexes $\mathcal{K}(0,i)$
associated with the fragment in example \ref{exa: ALP - G.H. a F.J.H. cc.1-2}.
\end{example}
\begin{figure}[btp]
\begin{centering}
\includegraphics[scale=0.9]{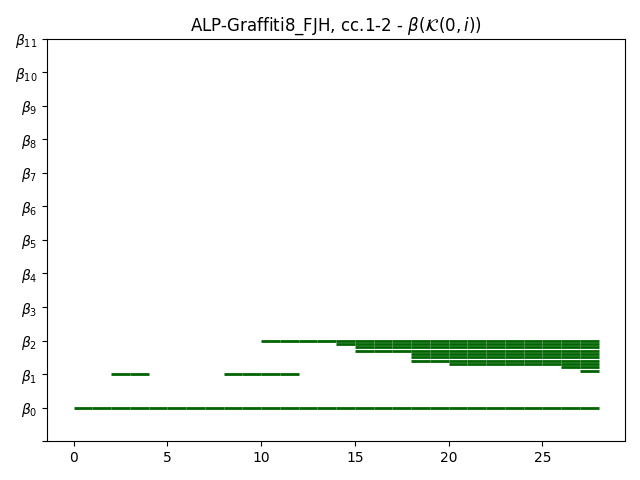}
\par\end{centering}

\protect\caption{{\label{fig:Barcodes K(0,i)  ALP-G.H. a F.J.H. cc1-2}Barcode
plot for Betti numbers of complexes $\mathcal{K}(0,i)$, for $0\leq i\leq27$,
from the $28\text{\textendash}$event fragment associated with Armando
Luna }- \emph{G.H. to F.J.H.}, mm. 1-2 (see figure \ref{fig: ALP - G.H. a F.J.H. cc. 1-2}
and example \ref{exa: ALP - G.H. a F.J.H. cc.1-2}).}

\end{figure}

\subparagraph{Simplicial complexes of cumulative events of radius $r$. \protect \\
}

\begin{figure}[btp]
\begin{centering}
\includegraphics[scale=0.9]{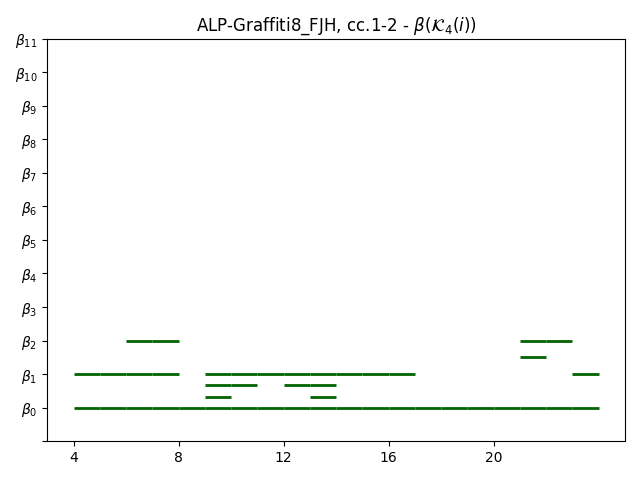}
\par\end{centering}

\protect\caption{{\label{fig:Barcodes K_4(i)  ALP-G.H. a F.J.H. cc1-2}Barcode
plot for Betti numbers of complexes for intervals of events of radius
$4$, $\mathcal{K}_{4}(i)$ for $4\leq i\leq24$, from fragment associated
with Armando Luna }- \emph{G.H. to F.J.H.}, mm. 1-2 (see figure \ref{fig: ALP - G.H. a F.J.H. cc. 1-2}
and example \ref{exa: ALP - G.H. a F.J.H. cc.1-2}).}

\end{figure}

As a special case of the above, given an integer $r\in\{0,\lfloor\frac{N}{2}\rfloor\}$
and a sequence of chords $\mathcal{A}(\mathcal{M})=\{a_{0},\,...\,,a_{N}\}$,
we consider the complex $\mathcal{K}_{r}(i)=\mathcal{K}(i-r,i+r)$,
corresponding to the interval 
\[
[a_{i-r},a_{i+r}]=\{a_{i-r},a_{i-r+1},\,...\,,a_{i-1},a_{i},a_{i+1},\,...\,,a_{i+r-1},a_{i+r}\}\,\text{.}
\]
We call the complex $\mathcal{K}_{r}(i)$ the \textbf{simplicial complex
of events of radius $\boldsymbol{r}$ }around $\mathbf{e}_{i}$ in
$\mathcal{M}$. These complexes give us local information about harmonic
sequences. Note that the resulting simplicial complexes for a fixed radius $r$ are not contained into one another. Thus, this sequence of simplicial complexes does not define a filtration of the complete complex $\mathcal K (0,N)$, but only a cover of it. As a consequence of this, we cannot strictly speak of persistent homology, and so the corresponding barcodes depict the values of Betti numbers of a sequence of simplicial complexes, without representing persitent homology cycles. Focusing on the subsequence of cumulative events of varying radii around a fixed event $e_{i_0}$,  
\[
\mathcal{K}_{0}(i_0) , \mathcal{K}_{1}(i_0)  , \mathcal{K}_{2}(i_0) , \mathcal{K}_{3}(i_0) , \, ... \,  \, ,
\]
we obtain an actual filtration. As an example, from the score in figure \ref{fig: ALP - G.H. a F.J.H. cc. 1-2}
we get the barcode plot shown in figure \ref{fig:Barcodes K_4(i)  ALP-G.H. a F.J.H. cc1-2},showing the Betti numbers of complexes of cumulative events of radius
$4$, $\mathcal{K}_{4}(i)$.

In a subsequent work we will focus on studying the results of calculating
these homological descriptors for all possible radii $r$. Note that
all simplicial complexes $\mathcal{K}(i,j)$ (including $\mathcal{K}_{r}(i)$)
contain information about pitch classes common to chords as well as
the number of pitch classes that constitue them. However, they do
not capture intervals in chords, which is a crucial stylistic feature.
In order to catch intervals in building simplicial complexes from
chords in a music fragment, we propose the construction described
in the next subsection.

\subsubsection{\label{sub:Simplicial-complexes-of chord by pitch and interval}Simplicial
complexes of cumulative events by pitch and interval. }

 We also consider another family of simplicial complexes on subsets
of pitch classes in the twelve tone equally tempered system. We associate
to a given chord $a$ with normal form vector $(\overline{x_{0}},\overline{x_{1}},\,...\,,\overline{x_{n}})$,
the simplicial complex whose simplices are $\sigma_{i}=\{\overline{x_{i}},\overline{x_{i}+1},\overline{x_{i}+2},\,...\,,\overline{x_{i+1}}\},\,i\in\{0,\,...\,,n-1\}$,
together with their faces, where pitch classes $\overline{x_{i}}\in\mathbb{Z}_{12}$
are always assumed to be represented by the smallest possible nonnegative
integer. Given a sequence of chords $\mathcal{A}=\{a_{0},a_{1},\,...\,,a_{N}\}$,
we denote the simplicial complex associated in this form to chord
$a_{i}$ by $\tilde{\mathcal{K}}(i,i)$, and proceed to define $\tilde{\mathcal{K}}(i,j)$
for $0\leq i\leq j\leq N$ similarly to how we defined $\mathcal{K}(i,j)$
(see the previous section). Thus, from this construction we get
a sequence of simplicial complexes representing the harmonic subsequences
of $\mathcal{A}$, upon which we can run a homology analysis.

To illustrate this construction, take again chord $(\bar{5},\bar{9},\bar{0})$
for example. With this construction we get a simplicial complex made
up from simplices $\{\bar{5},\bar{6},\bar{7},\bar{8},\bar{9}\}$ and
$\{\bar{9},\overline{10},\overline{11},\bar{0}\}$, which represent
the major third interval between pitch classes $\bar{5}$ and $\bar{9}$,
and the minor third between $\bar{9}$ and $\bar{0}$, respectively,
together with all their faces. Considering this chord is the beginning
of fragment from example \ref{exa: ALP - G.H. a F.J.H. cc.1-2}, we
get $\tilde{\mathcal{K}}(0,0)$ for that score (see figure \ref{fig: Simplicial complex by pitch and interval for chord (5,9,0)}).

\begin{figure}[btp]
\begin{centering}
\includegraphics[scale=0.5]{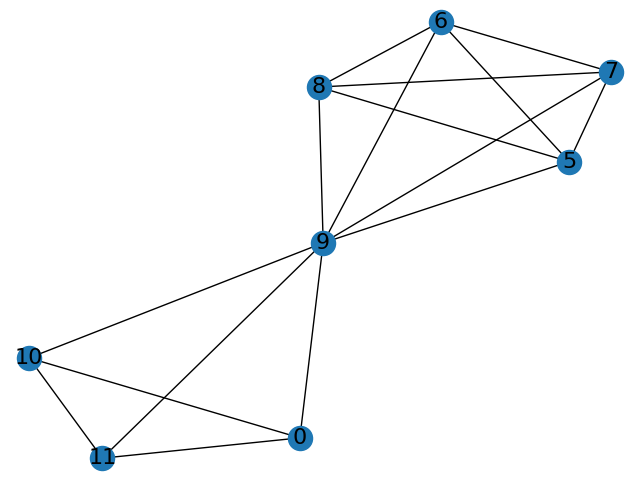}
\par\end{centering}

\protect\caption{\label{fig: Simplicial complex by pitch and interval for chord (5,9,0)}{Simplicial
complex of chord $(\bar{5},\bar{9},\bar{0})$ by pitch class and interval.} }

\end{figure}

\section{Results. \label{sec:Results.}}

To test our proposal as a way of describing and comparing musical
data, we analyzed the persistence diagrams corresponding to the different
data mappings associated with four fragments taken from different musical
examples sharing some common stylistic elements. Such fragments belong
to the following pieces: another one of Armando Luna's \emph{Graffiti}
(see first paragraph of example \ref{exa: ALP - G.H. a F.J.H. cc.1-2}),
the one dedicated to J.S. Bach (\emph{G.H. to J.S.B.}), and three
actual pieces by J.S. Bach: \emph{Brandenburg Concertos nos. 1-3 BWV
1046-1048} (\emph{B.C. 1-3}). For this test we consider only the first
four measures of each piece. We compare the harmonic data (with and
without temporal data) contained in these examples by using the bottleneck
distance calculated between their $H_{0}\text{\textendash}$ and $H_{1}\text{\textendash}$diagrams
for all six data mappings (see tables \ref{tab: H_0 Bottleneck-distances ALP-G.H. a J.S.B. vs 3 Brandenburg}
and \ref{tab: H_1 Bottleneck-distances ALP-G.H. a J.S.B. vs 3 Brandenburg}).
Diagrams for homology in dimensions greater than $0$ ($H_{1},H_{2},\,...\,$)
appear jointly only for some mappings and some pairs of samples\footnote{$H_{0}\text{\textendash}$diagrams are necessarily non empty, as a
finite data set is a bounded set in some $\mathbb{R}^{n}$, and so
in the successive construction of Vietoris-Rips complexes eventually
at least one connected component is always persistent.}, so it is possible to compare these diagrams in only a few cases.
For the particular scores analyzed, we get the possibility of comparing
their $H_{1}\text{\textendash}$diagrams for a few mappings (no $H_{2}\text{\textendash}$
or higher dimension diagrams were jointly generated for any pair of
these examples; refer to table \ref{tab: H_1 Bottleneck-distances ALP-G.H. a J.S.B. vs 3 Brandenburg}).
We do not reproduce here all persistence and barcode diagrams associated
with our examples, but present the dendrograms that show the comparisons
among $H_{0}\text{\textendash}$ and $H_{1}\text{\textendash}$diagrams.

\begin{figure}[btp]
\begin{centering}
\includegraphics[scale=0.5]{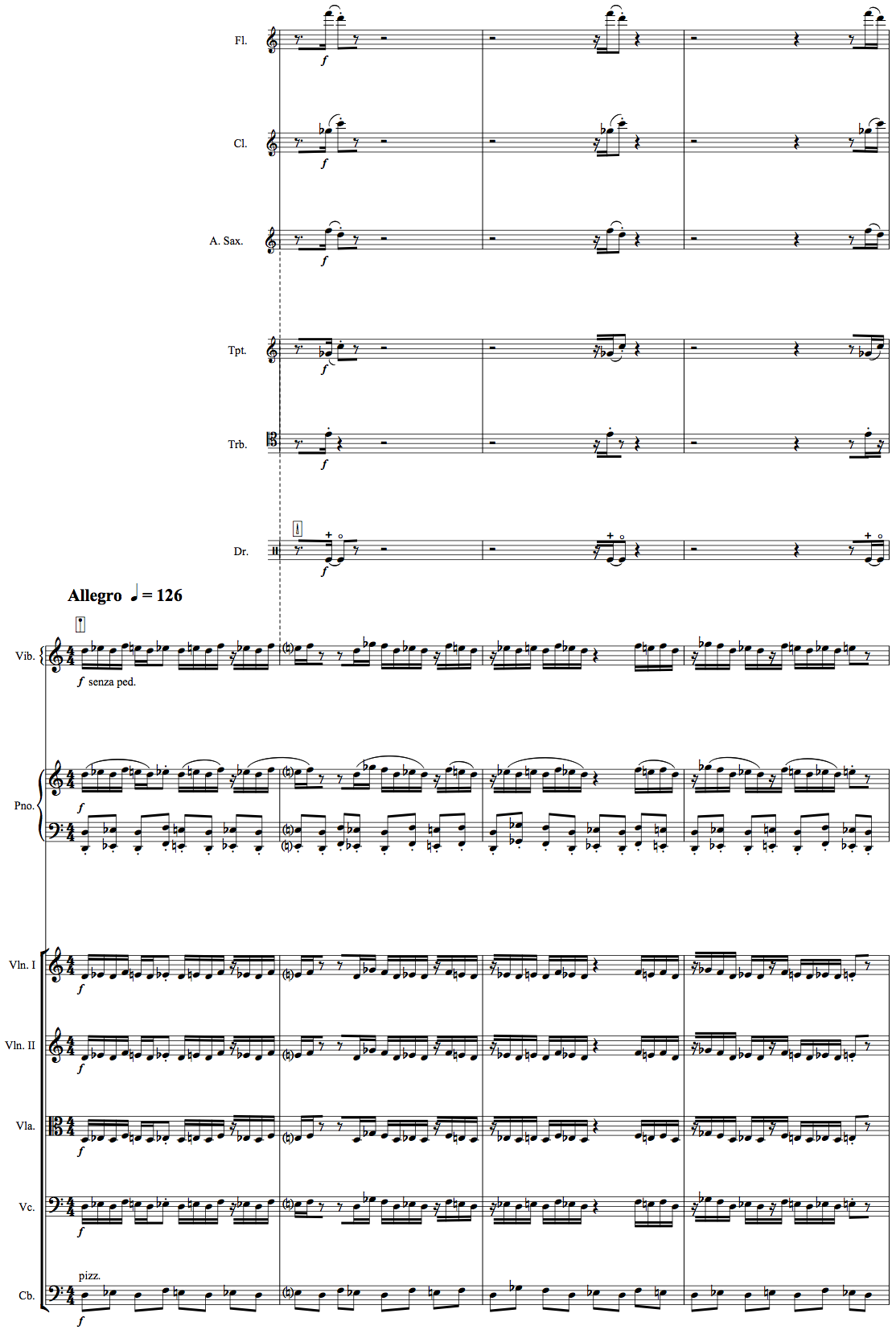}
\par\end{centering}

\protect\caption{{\label{fig: ALP-Graffiti1_JSB_cc1-4}Measures 1-4
of Armando Luna's }\emph{Graffiti Hommage a J.S. Bach}.}

\end{figure}

In the following, we go through the dendrograms (hierarchical clustering
plots) depicting bottleneck distances between $H_{0}\text{\textendash}$persistence
diagrams from all six data mappings for fragments corresponding to
the first four measures of Luna's \emph{Graffiti Hommage to J.S. Bach}
(see figure \ref{fig: ALP-Graffiti1_JSB_cc1-4}) and J.S. Bach's first
three \emph{Brandenburg Concertos}. For the sake of completeness,
we also present dendrograms for bottleneck distances between $H_{1}\text{\textendash}$persistence
diagrams from data mappings I, II, IV, V, and VI (the mappings for
which there are $H_{1}\text{\textendash}$diagrams for all of our
chosen samples; see table \ref{tab: H_1 Bottleneck-distances ALP-G.H. a J.S.B. vs 3 Brandenburg}).
Bottleneck distance values are represented on the vertical axis of
these plots, while numbers above the samples' identification on the
horizontal axis are only labels.

\begin{center}
\begin{table}[!tph]
\begin{centering}

{\tiny{}}%
\resizebox{\textwidth}{!}{\begin{tabular}{|c|c|c|c|c|c|c|c|c|c|c|c|c|c|c|c|c|c|c|}
\hline 
 & \multicolumn{6}{c|}{{\tiny{}\cellcolor{lightgray}B.C.1-BWV1046}} & \multicolumn{6}{c|}{{\tiny{}\cellcolor{lightgray}B.C.2-BWV1047}} & \multicolumn{6}{c|}{{\tiny{}\cellcolor{lightgray}B.C.3-BWV1048}}\tabularnewline
\hline 
\hline 
{\tiny{}Data mapping:} & {\tiny{}\cellcolor{purplepizzazz}I} & {\tiny{}\cellcolor{cyan}II} & {\tiny{}\cellcolor{yellow}III} & {\tiny{}\cellcolor{orange}IV} & {\tiny{}\cellcolor{green}V} & {\tiny{}\cellcolor{blue-magenta}VI} & {\tiny{}\cellcolor{purplepizzazz}I} & {\tiny{}\cellcolor{cyan}II} & {\tiny{}\cellcolor{yellow}III} & {\tiny{}\cellcolor{orange}IV} & {\tiny{}\cellcolor{green}V} & {\tiny{}\cellcolor{blue-magenta}VI} & {\tiny{}\cellcolor{purplepizzazz}I} & {\tiny{}\cellcolor{cyan}II} & {\tiny{}\cellcolor{yellow}III} & {\tiny{}\cellcolor{orange}IV} & {\tiny{}\cellcolor{green}V} & {\tiny{}\cellcolor{blue-magenta}VI}\tabularnewline
\hline 
{\tiny{}\cellcolor{lightgray}ALP-G.H. to J.S.B.} & {\tiny{}\cellcolor{purplepizzazz}6.08} & {\tiny{}\cellcolor{cyan}1.08} & {\tiny{}\cellcolor{yellow}1} & {\tiny{}\cellcolor{orange}0.5} & {\tiny{}\cellcolor{green}6.2} & {\tiny{}\cellcolor{blue-magenta}1.5} & {\tiny{}\cellcolor{purplepizzazz}6.62} & {\tiny{}\cellcolor{cyan}1.45} & {\tiny{}\cellcolor{yellow}0.71} & {\tiny{}\cellcolor{orange}0.5} & {\tiny{}\cellcolor{green}6.4} & {\tiny{}\cellcolor{blue-magenta}1.58} & {\tiny{}\cellcolor{purplepizzazz}8.72} & {\tiny{}\cellcolor{cyan}1.95} & {\tiny{}\cellcolor{yellow}1.03} & {\tiny{}\cellcolor{orange}0.5} & {\tiny{}\cellcolor{green}8.87} & {\tiny{}\cellcolor{blue-magenta}1.5}\tabularnewline
\hline 
{\tiny{}\cellcolor{lightgray}B.C.1-BWV1046} & \multicolumn{6}{c|}{} & {\tiny{}\cellcolor{purplepizzazz}5.07} & {\tiny{}\cellcolor{cyan}0.72} & {\tiny{}\cellcolor{yellow}0.71} & {\tiny{}\cellcolor{orange}0.5} & {\tiny{}\cellcolor{green}4.85} & {\tiny{}\cellcolor{blue-magenta}1.12} & {\tiny{}\cellcolor{purplepizzazz}7.39} & {\tiny{}\cellcolor{cyan}0.86} & {\tiny{}\cellcolor{yellow}0.82} & {\tiny{}\cellcolor{orange}0.5} & {\tiny{}\cellcolor{green}6.08} & {\tiny{}\cellcolor{blue-magenta}1.23}\tabularnewline
\cline{1-1} \cline{8-19} 
{\tiny{}\cellcolor{lightgray}B.C.2-BWV1047} & \multicolumn{6}{c}{} & \multicolumn{6}{c|}{} & {\tiny{}\cellcolor{purplepizzazz}6.48} & {\tiny{}\cellcolor{cyan}0.5} & {\tiny{}\cellcolor{yellow}0.59} & {\tiny{}\cellcolor{orange}0.5} & {\tiny{}\cellcolor{green}6.42} & {\tiny{}\cellcolor{blue-magenta}1.22}\tabularnewline
\cline{1-1} \cline{14-19} 
\end{tabular}}
\par\end{centering}{\tiny \par}

\begin{centering}
{\footnotesize{}\protect\caption{{\scriptsize{}\label{tab: H_0 Bottleneck-distances ALP-G.H. a J.S.B. vs 3 Brandenburg}Bottleneck
distances (rounded up to two decimal positions) between $H_{0}\text{\textendash}$diagrams
from six data mappings for mm. 1\textendash 4 of Luna's }\emph{\scriptsize{}Graffiti
Hommage to J.S.B.}{\scriptsize{} and J.S. Bach's }\emph{\scriptsize{}Brandenburg
Concertos BWV 1046-1048}.}
}
\par\end{centering}{\footnotesize \par}

\end{table}

\par\end{center}

\begin{center}

\par\end{center}

\begin{center}
\begin{table}[!tph]
\begin{centering}
{\tiny{}}%
\resizebox{\textwidth}{!}{\begin{tabular}{|c|c|c|c|c|c|c|c|c|c|c|c|c|c|c|c|c|c|c|}
\hline 
 & \multicolumn{6}{c|}{{\tiny{}\cellcolor{lightgray}B.C.1-BWV1046}} & \multicolumn{6}{c|}{{\tiny{}\cellcolor{lightgray}B.C.2-BWV1047}} & \multicolumn{6}{c|}{{\tiny{}\cellcolor{lightgray}B.C.3-BWV1048}}\tabularnewline
\hline 
\hline 
{\tiny{}Data mapping:} & {\tiny{}\cellcolor{purplepizzazz}I} & {\tiny{}\cellcolor{cyan}II} & {\tiny{}\cellcolor{yellow}III} & {\tiny{}\cellcolor{orange}IV} & {\tiny{}\cellcolor{green}V} & {\tiny{}\cellcolor{blue-magenta}VI} & {\tiny{}\cellcolor{purplepizzazz}I} & {\tiny{}\cellcolor{cyan}II} & {\tiny{}\cellcolor{yellow}III} & {\tiny{}\cellcolor{orange}IV} & {\tiny{}\cellcolor{green}V} & {\tiny{}\cellcolor{blue-magenta}VI} & {\tiny{}\cellcolor{purplepizzazz}I} & {\tiny{}\cellcolor{cyan}II} & {\tiny{}\cellcolor{yellow}III} & {\tiny{}\cellcolor{orange}IV} & {\tiny{}\cellcolor{green}V} & {\tiny{}\cellcolor{blue-magenta}VI}\tabularnewline
\hline 
{\tiny{}\cellcolor{lightgray}ALP-G.H. to J.S.B.} & {\tiny{}\cellcolor{purplepizzazz}0.64} & {\tiny{}\cellcolor{cyan}0.55} &  & {\tiny{}\cellcolor{orange}0.21} & {\tiny{}\cellcolor{green}0.55} & {\tiny{}\cellcolor{blue-magenta}0.56} & {\tiny{}\cellcolor{purplepizzazz}0.76} & {\tiny{}\cellcolor{cyan}0.28} &  & {\tiny{}\cellcolor{orange}0.21} & {\tiny{}\cellcolor{green}0.21} & {\tiny{}\cellcolor{blue-magenta}0.42} & {\tiny{}\cellcolor{purplepizzazz}0.39} & {\tiny{}\cellcolor{cyan}0.27} &  & {\tiny{}\cellcolor{orange}0.21} & {\tiny{}\cellcolor{green}0.31} & {\tiny{}\cellcolor{blue-magenta}0.5}\tabularnewline
\cline{1-3} \cline{5-6} \cline{8-9} \cline{11-15} \cline{17-19} 
{\tiny{}\cellcolor{lightgray}B.C.1-BWV1046} & \multicolumn{6}{c|}{} & {\tiny{}\cellcolor{purplepizzazz}0.65} & {\tiny{}\cellcolor{cyan}0.55} &  & {\tiny{}\cellcolor{orange}0} & {\tiny{}\cellcolor{green}0.55} & {\tiny{}\cellcolor{blue-magenta}0.56} & {\tiny{}\cellcolor{purplepizzazz}6.44} & {\tiny{}\cellcolor{cyan}0.55} &  & {\tiny{}\cellcolor{orange}0.21} & {\tiny{}\cellcolor{green}0.55} & {\tiny{}\cellcolor{blue-magenta}0.59}\tabularnewline
\cline{1-1} \cline{8-9} \cline{11-15} \cline{17-19} 
{\tiny{}\cellcolor{lightgray}B.C.2-BWV1047} & \multicolumn{6}{c}{} & \multicolumn{6}{c|}{} & {\tiny{}\cellcolor{purplepizzazz}0.76} & {\tiny{}\cellcolor{cyan}0.21} &  & {\tiny{}\cellcolor{orange}0.21} & {\tiny{}\cellcolor{green}0.31} & {\tiny{}\cellcolor{blue-magenta}0.59}\tabularnewline
\cline{1-1} \cline{14-15} \cline{17-19} 
\end{tabular}}
\par\end{centering}{\tiny \par}

\centering{}{\footnotesize{}\protect\caption{{\scriptsize{}\label{tab: H_1 Bottleneck-distances ALP-G.H. a J.S.B. vs 3 Brandenburg}Bottleneck
distances (rounded up to two decimal positions) between $H_{1}\text{\textendash}$diagrams
from six data mappings for mm. 1\textendash 4 of Luna's }\emph{\scriptsize{}Graffiti
Hommage to J.S.B.}{\scriptsize{} and J.S. Bach's }\emph{\scriptsize{}Brandenburg
Concertos BWV 1046-1048}.}
}{\footnotesize \par}
\end{table}

\par\end{center}

In order to contrast our approach with the viewpoint of traditional
harmony, we synthesize chord progressions, in traditional tonal nomenclature,
for \emph{Brandenburg Concertos 1-3} in table \ref{tab: Brandenburgs 1-3 acordes}
(we complete the table with a tonal interpretation of Luna's piece).
When interpreting these progressions against our measurements, one
must keep in mind that events may not always reflect this chords literally,
as they may include harmonic ornaments such as passing notes, auxiliary
notes, retardations, etc., which are left out of the traditional tonal
notation.

\begin{figure}[btp]
\begin{centering}
\includegraphics[scale=0.37]{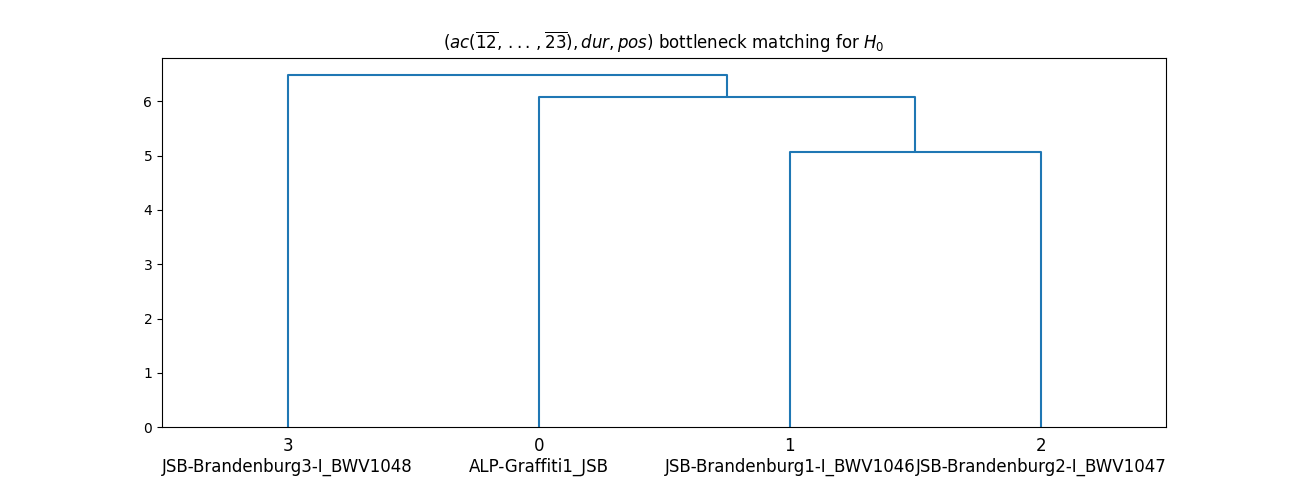}
\par\end{centering}

\protect\caption{{\scriptsize{}\label{fig:Dendrogram-for-bottleneck_H0 cat 1}Dendrogram for bottleneck distances between $H_{0}\text{\textendash}$diagrams
for }\textbf{\scriptsize{}data mapping I}{\scriptsize{} from mm. 1\textendash 4
of Luna's }\emph{\scriptsize{}Graffiti Hommage to J.S.B.}{\scriptsize{}
and J.S. Bach's }\emph{\scriptsize{}Brandenburg Concertos BWV 1046-1048}.}
\end{figure}

\begin{center}
\begin{table}[!tph]
\begin{centering}
{\tiny{}}%
\resizebox{\textwidth}{!}{\begin{tabular}{|c|c|c|c|c|}
\hline 
 & {\footnotesize{}\cellcolor{lightgray}m. 1} & {\footnotesize{}\cellcolor{lightgray}m. 2} & {\footnotesize{}\cellcolor{lightgray}m. 3} & {\footnotesize{}\cellcolor{lightgray}m. 4}\tabularnewline
\hline 
\hline 
{\tiny{}\cellcolor{lightgray}B.C.1-BWV1046 (FMaj)} & I-VI-II-V7 & I-V-VI7-IIM7 & V-VI7 & I-V7-I\tabularnewline
\hline 
{\tiny{}\cellcolor{lightgray}B.C.2-BWV1047 (FMaj)} & I & I & I-V7-I-V & III7-V7-I-V7\tabularnewline
\hline 
{\tiny{}\cellcolor{lightgray}B.C.2-BWV1048 (GMaj)} & I (upbeat) & I & I-VI-IIM7-V-I\musFig{6 4} & I\musFig{6 4}-V-V7\tabularnewline
\hline 
{\tiny{}\cellcolor{lightgray}ALP-G.H. to J.S.B. (``D min/Maj'')} & ``I'' & ``Im/M''-''I'' & ``Im/M''-''I'' & ``Im/M''-''I''\tabularnewline
\hline 
\end{tabular}}
\par\end{centering}{\tiny \par}

\centering{}{\footnotesize{}\protect\caption{{\scriptsize{}\label{tab: Brandenburgs 1-3 acordes}Chord progressions
in mm. 1\textendash 4 of Luna's }\emph{\scriptsize{}Graffiti Hommage
to J.S.B.}{\scriptsize{} and J.S. Bach's }\emph{\scriptsize{}Brandenburg
Concertos BWV 1046-1048}.}
}{\footnotesize \par}
\end{table}

\par\end{center}
\begin{itemize}
\item In the plot shown in figure  \ref{fig:Dendrogram-for-bottleneck_H0 cat 1} we see the distances between the $H_{0}\text{\textendash}$persistence
diagrams describing the clouds of data points resulting from mapping
chords with their rhythms and onsets into $\mathbb{R}^{14}$,
as tuples which include pitch information in their first twelve entries
as integers from $12$ to $23$ (instead of $0$ through $12$). We
observe that according to this embedding of events into $\mathbb{R}^{14}$,
the first four bars of the first movements of \emph{Brandenburg Concertos
1 }and \emph{2 }are the closest in shape among all four samples. This
is compatible with the traditional tonal analysis of these pieces,
since the first movements of \emph{Brandeburg Concertos 1 }and \emph{2}
are both in the key of F Major (as we said before, this data mapping
is sensitive to transpositions), and have a very similar distribution
of rhythms (mostly eights and sixteenths). Nevertheless, these two
movements differ notably in their harmonic march, which will be discussed
below. On the other hand, Luna's \emph{Graffiti }and \emph{Brandenburg
Concerto no. 3 }appear as the next pair of closest samples, although
\emph{B.C.3 }is the furthest. This may be explained by the fact that
\emph{Brandenburg 3 }is written in the G Major key, while Luna's \emph{Graffiti
}includes many instances of D minor (the relative minor of F) and
D Major (the dominant of G).
\end{itemize}

\begin{figure}[btp]
\begin{centering}
\includegraphics[scale=0.37]{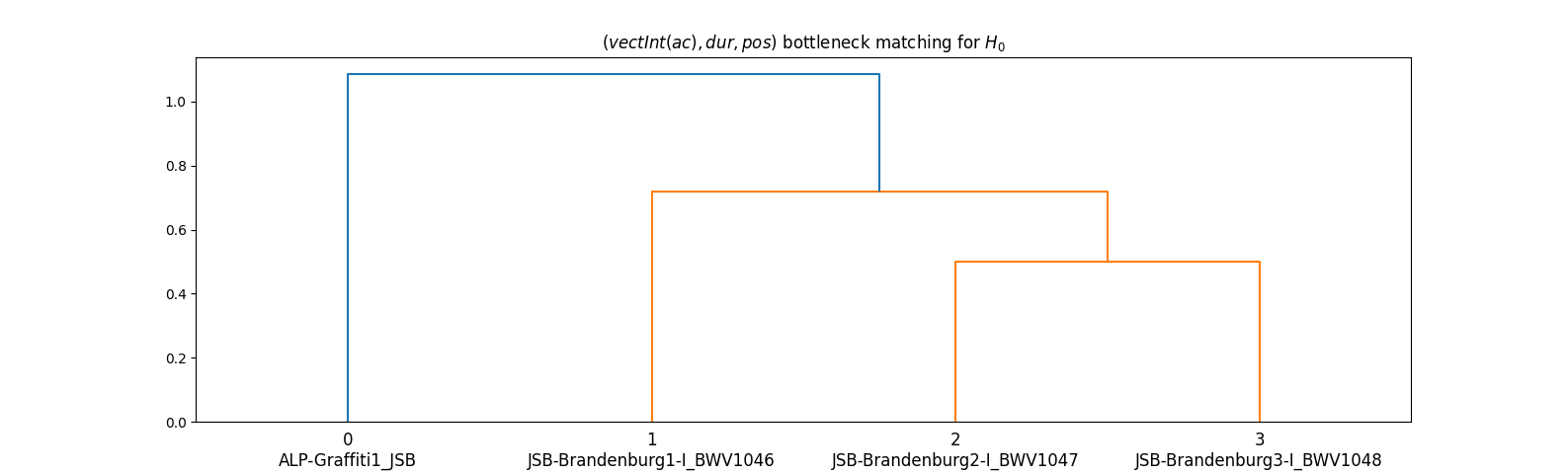}
\par\end{centering}

\protect\caption{{\scriptsize{}\label{fig:Dendrogram-for-bottleneck_H0 cat 2}Dendrogram for bottleneck distances between $H_{0}\text{\textendash}$diagrams
for }\textbf{\scriptsize{}data mapping II}{\scriptsize{} from mm.
1\textendash 4 of Luna's }\emph{\scriptsize{}Graffiti Hommage to J.S.B.}{\scriptsize{}
and J.S. Bach's }\emph{\scriptsize{}Brandenburg Concertos BWV 1046-1048}.}
\end{figure}

\begin{itemize}
\item In figure \ref{fig:Dendrogram-for-bottleneck_H0 cat 2}, in which interval vectors instead of normal form
vectors are considered (together with durations and onsets
of events), we observe that \emph{Brandenburg Concertos 2} and \emph{3
}are the most similar samples, with \emph{Brandenburg 1 }being the
next one closest to them, and Luna's \emph{Graffiti} as substantially
dissimilar to the rest of the samples. This may be explained by the
auxiliary notes present over the triads in these pieces, which are
mostly fourths, minor sevenths, and major seconds in Bach, and minor
seconds in Luna. Thus, it makes sense that, beyond the tonality of
the fragments, harmonic, and rhythmic content in both authors is fairly
distinguishable.
\end{itemize}

\begin{figure}[btp]
\begin{centering}
\includegraphics[scale=0.37]{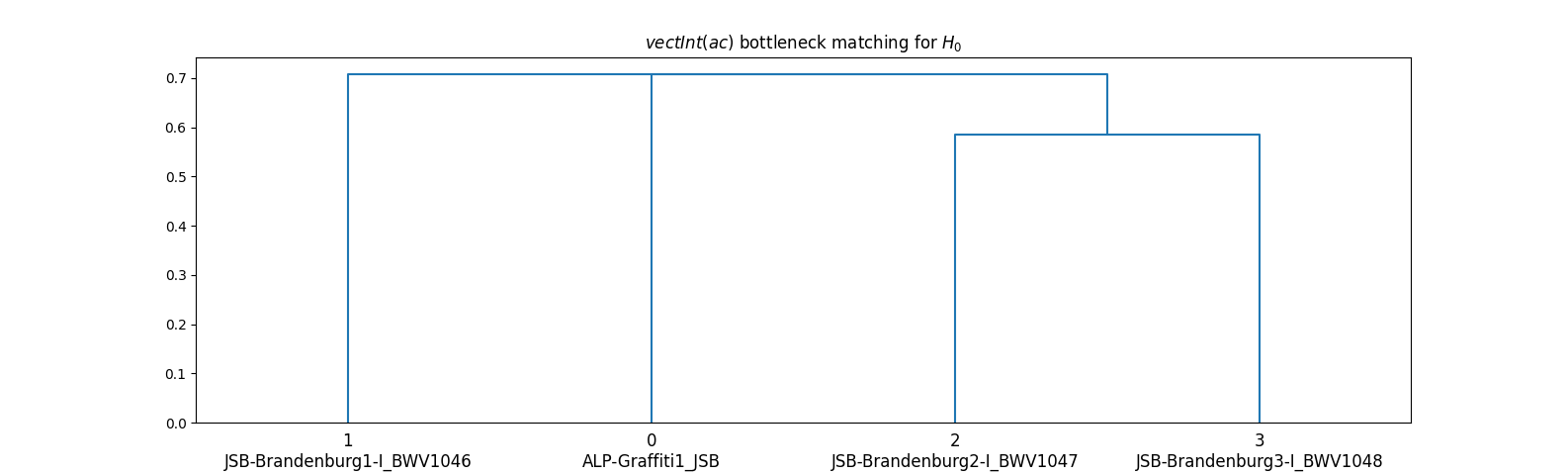}
\par\end{centering}

\protect\caption{{\scriptsize{}\label{fig:Dendrogram-for-bottleneck_H0 cat 3}Dendrogram
for bottleneck distances between $H_{0}\text{\textendash}$diagrams
for }\textbf{\scriptsize{}data mapping III}{\scriptsize{} from mm.
1\textendash 4 of Luna's }\emph{\scriptsize{}Graffiti Hommage to J.S.B.}{\scriptsize{}
and J.S. Bach's }\emph{\scriptsize{}Brandenburg Concertos BWV 1046-1048}.}
\end{figure}

\begin{itemize}
\item The third data mapping, formed by interval vectors only, produces
the distances shown in figure \ref{fig:Dendrogram-for-bottleneck_H0 cat 3}
among $H_{0}\text{\textendash}$diagrams of the analyzed fragments.
We can conclude from it almost the same relations as from the plot
for mapping II, except that in this case \emph{Brandenburg Concerto
1} and \emph{Graffiti Hommage to J.S.B. }are equidistant to the cluster
formed by \emph{Brandenburg Concertos 1} and \emph{3}. This tells
us that vertical intervalic relations are similar in \emph{B.C. 1}
and \emph{G.H. to J.S.B.}, showing similarity in their harmonic styles.
\end{itemize}

\begin{figure}[btp]
\begin{centering}
\includegraphics[scale=0.37]{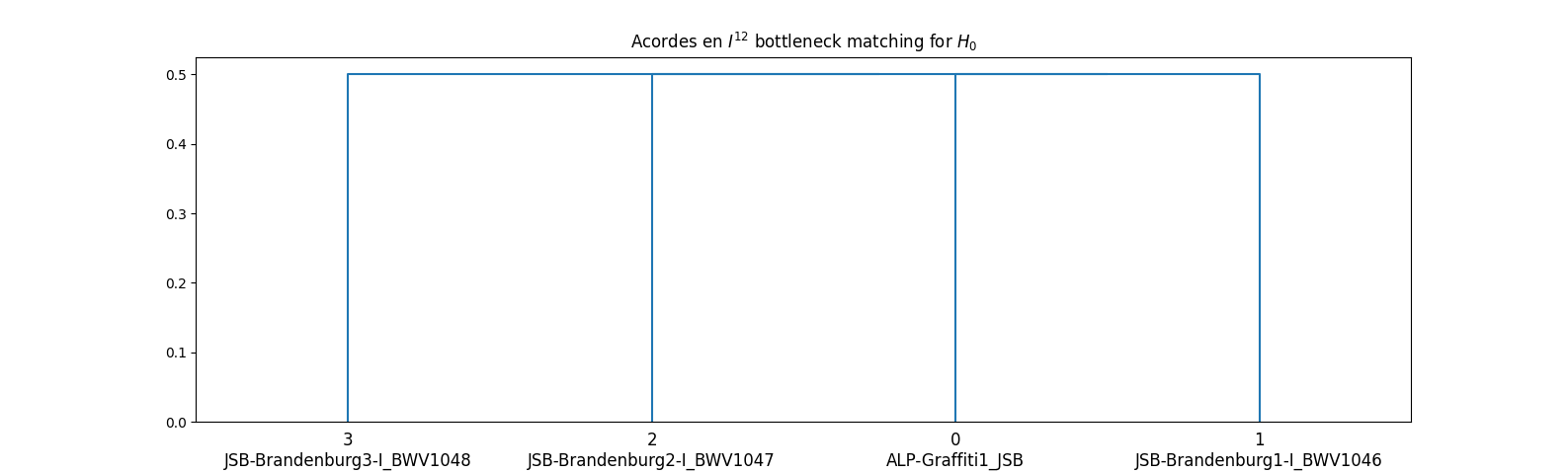}
\par\end{centering}

\protect\caption{{\scriptsize{}\label{fig:Dendrogram-for-bottleneck_H0 cat 4}Dendrogram for bottleneck distances between $H_{0}\text{\textendash}$diagrams
for }\textbf{\scriptsize{}data mapping IV}{\scriptsize{} from mm.
1\textendash 4 of Luna's }\emph{\scriptsize{}Graffiti Hommage to J.S.B.}{\scriptsize{}
and J.S. Bach's }\emph{\scriptsize{}Brandenburg Concertos BWV 1046-1048}.}
\end{figure}

\begin{figure}[btp]
\begin{centering}
\includegraphics[scale=0.37]{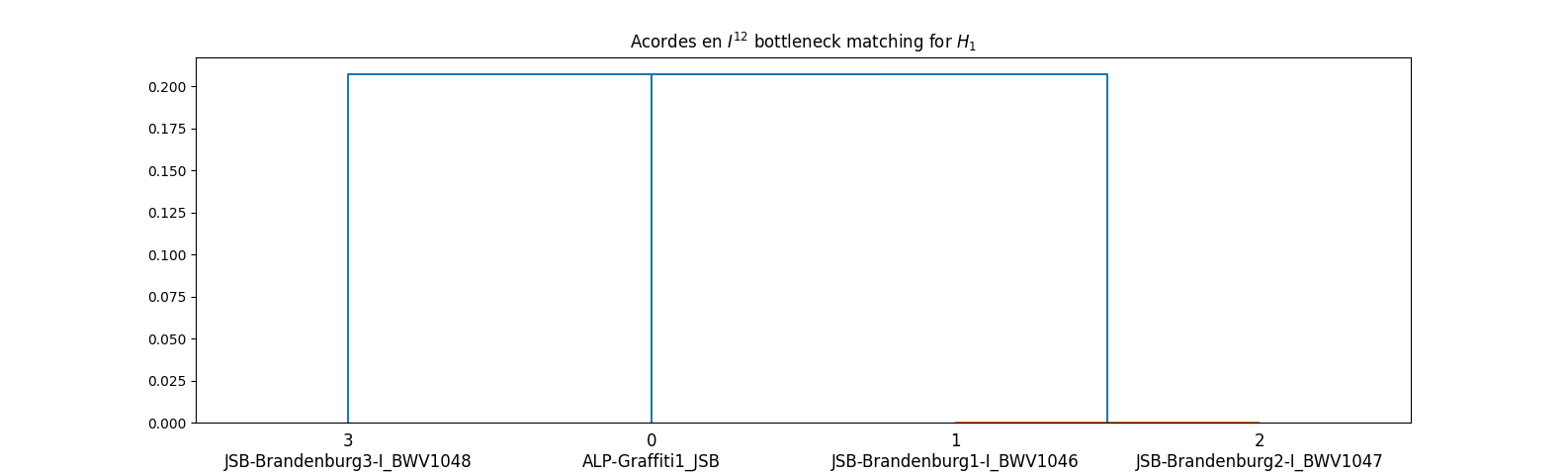}
\par\end{centering}

\protect\caption{{\scriptsize{}\label{fig:Dendrogram-for-bottleneck_H1 cat 4}Dendrogram
for bottleneck distances between $H_{1}\text{\textendash}$diagrams
for }\textbf{\scriptsize{}data mapping IV}{\scriptsize{} from mm.
1\textendash 4 of Luna's }\emph{\scriptsize{}Graffiti Hommage to J.S.B.}{\scriptsize{}
and J.S. Bach's }\emph{\scriptsize{}Brandenburg Concertos BWV 1046-1048}.}
\end{figure}

\begin{itemize}
\item The present samples are indistinguishable by comparing their $H_{0}\text{\textendash}$diagrams
associated with data mapping IV (see figure See figure \ref{fig:Dendrogram-for-bottleneck_H0 cat 4}). Nevertheless, this is not the case for $H_{1}\text{\textendash}$diagrams;
see figure \ref{fig:Dendrogram-for-bottleneck_H1 cat 4} ahead). This points
towards a similar structure among sets of vertical pitch sets.
\end{itemize}

\begin{figure}[btp]
\begin{centering}
\includegraphics[scale=0.37]{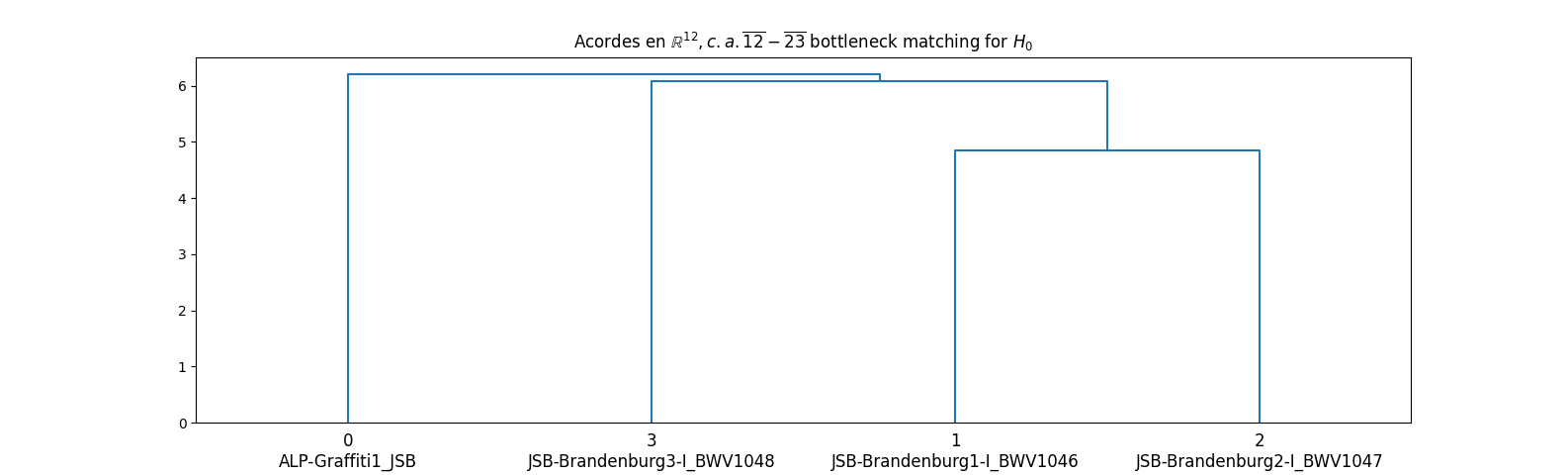}
\par\end{centering}

\protect\caption{{\scriptsize{}\label{fig:Dendrogram-for-bottleneck_H0 cat 5}Dendrogram for bottleneck distances between $H_{0}\text{\textendash}$diagrams
for }\textbf{\scriptsize{}data mapping V}{\scriptsize{} from mm. 1\textendash 4
of Luna's }\emph{\scriptsize{}Graffiti Hommage to J.S.B.}{\scriptsize{}
and J.S. Bach's }\emph{\scriptsize{}Brandenburg Concertos BWV 1046-1048}.}
\end{figure}

\begin{itemize}
\item As with diagrams coming from data mapping I, in data mapping V (see figure \ref{fig:Dendrogram-for-bottleneck_H0 cat 5}) the
three \emph{Brandenburgs }are closer among themselves than they are
to Luna's \emph{G.H. to J.S.B.}, being \emph{B.C. 1 }and\emph{ 2}
the most similar fragments. 
\end{itemize}

\begin{itemize}
\item Finally, for chords mapped as tuples in which de $i\text{\textendash}$th
coordinate represents the interval class in the chords' normal form
above pitch class $\bar{i}$, we get figure \ref{fig:Dendrogram-for-bottleneck ALP-JSB vs Brandenburgs 1-3 cc 1-4, cat 6}.
For this mapping, the clustering of the given samples coincides with
that of mappings II and III.
\end{itemize}
\begin{figure}[btp]
\begin{centering}
\includegraphics[scale=0.37]{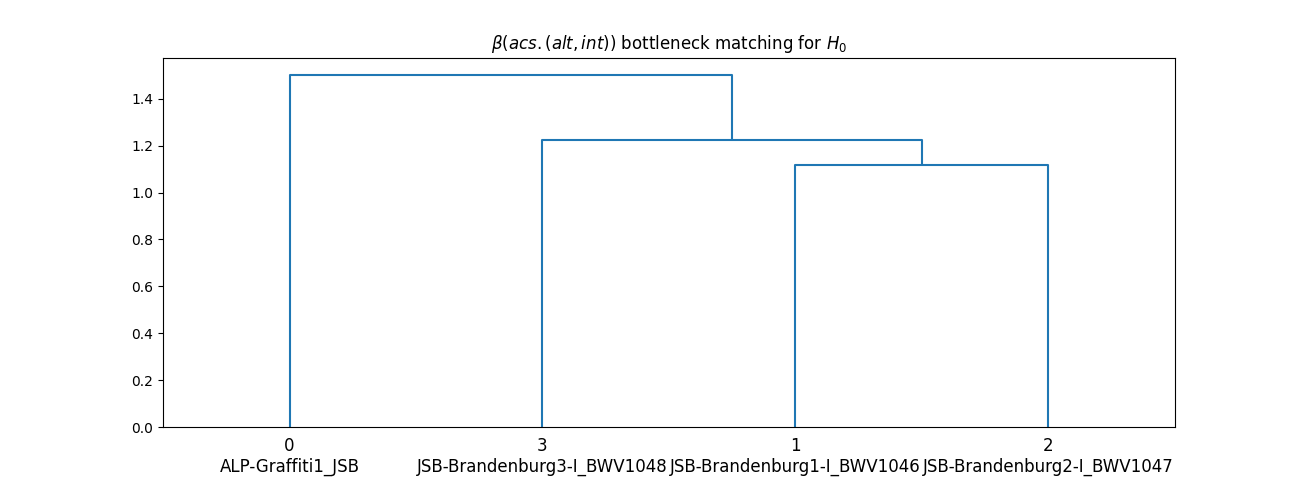}
\par\end{centering}

\protect\caption{{\scriptsize{}\label{fig:Dendrogram-for-bottleneck ALP-JSB vs Brandenburgs 1-3 cc 1-4, cat 6}Dendrogram
for bottleneck distances between $H_{0}\text{\textendash}$diagrams
for }\textbf{\scriptsize{}data mapping VI}{\scriptsize{} from mm.
1\textendash 4 of Luna's }\emph{\scriptsize{}Graffiti Hommage to J.S.B.}{\scriptsize{}
and J.S. Bach's }\emph{\scriptsize{}Brandenburg Concertos BWV 1046-1048}.}
\end{figure}

Now we present the dendrograms corresponding to $H_{1}\text{\textendash}$diagrams
from data mappings I, II, IV, V, and VI ($H_{1}\text{\textendash}$diagrams
do not exist for all of the samples for mapping III). From these we
may draw the following:
\begin{itemize}
\item For mapping I we get a kind of inversion in the ``closeness'' of
samples: for $H_{1}\text{\textendash}$diagrams, the closest samples
are \emph{G.H. to J.S.B.} and \emph{B.C. 3}, while \emph{B.C. 1 }and
\emph{2} generate clusters within a small distance. This suggests
some measure of similarity between data from a piece by Bach and Luna's
hommage.
\item Mapping II remains consistent in the closest pair of samples (though
not in the furthest one). 
\item Clustering in mapping IV is disambiguated, as the $H_{0}\text{\textendash}$diagrams
of all fragments were all equidistant for this mapping, while the
distance between $H_{1}\text{\textendash}$diagrams of \emph{B.C.
1} and \emph{2} is $0$.
\item $H_{1}\text{\textendash}$diagrams for mappings V and VI point towards
similarity between data from \emph{G.H. to J.S.B. }and \emph{B.C.
2}. Although from table \ref{tab: Interpretacion dendrograma ALP-G.H. a J.S.B. vs 3 Brandenburg}
the results obtained for this two mappings seem redundant, they are
not necessarily so for our analysis (let us not forget that mapping
V is sensitive to transpositions, while mapping VI is not; also mapping
VI incorporates intervals, while mapping V does not).
\end{itemize}
\begin{figure}[btp]
\begin{centering}
\includegraphics[scale=0.37]{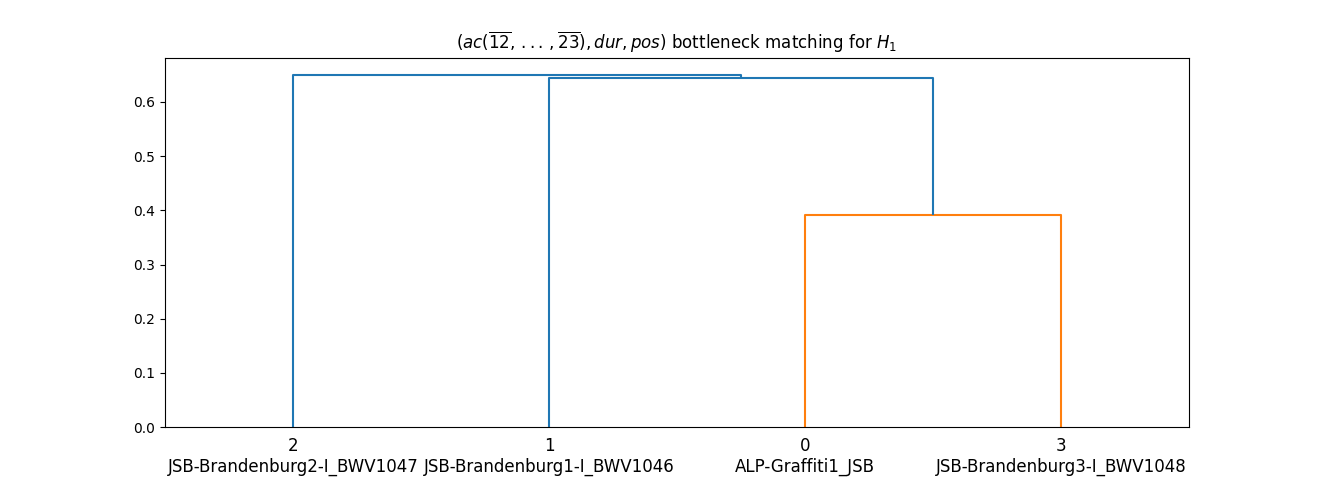}
\par\end{centering}

\protect\caption{{\scriptsize{}\label{fig:Dendrogram-for-bottleneck_H1 cat 1}Dendrogram
for bottleneck distances between $H_{1}\text{\textendash}$diagrams
for }\textbf{\scriptsize{}data mapping I}{\scriptsize{} from mm. 1\textendash 4
of Luna's }\emph{\scriptsize{}Graffiti Hommage to J.S.B.}{\scriptsize{}
and J.S. Bach's }\emph{\scriptsize{}Brandenburg Concertos BWV 1046-1048}.}
\end{figure}

\begin{figure}[btp]
\begin{centering}
\includegraphics[scale=0.37]{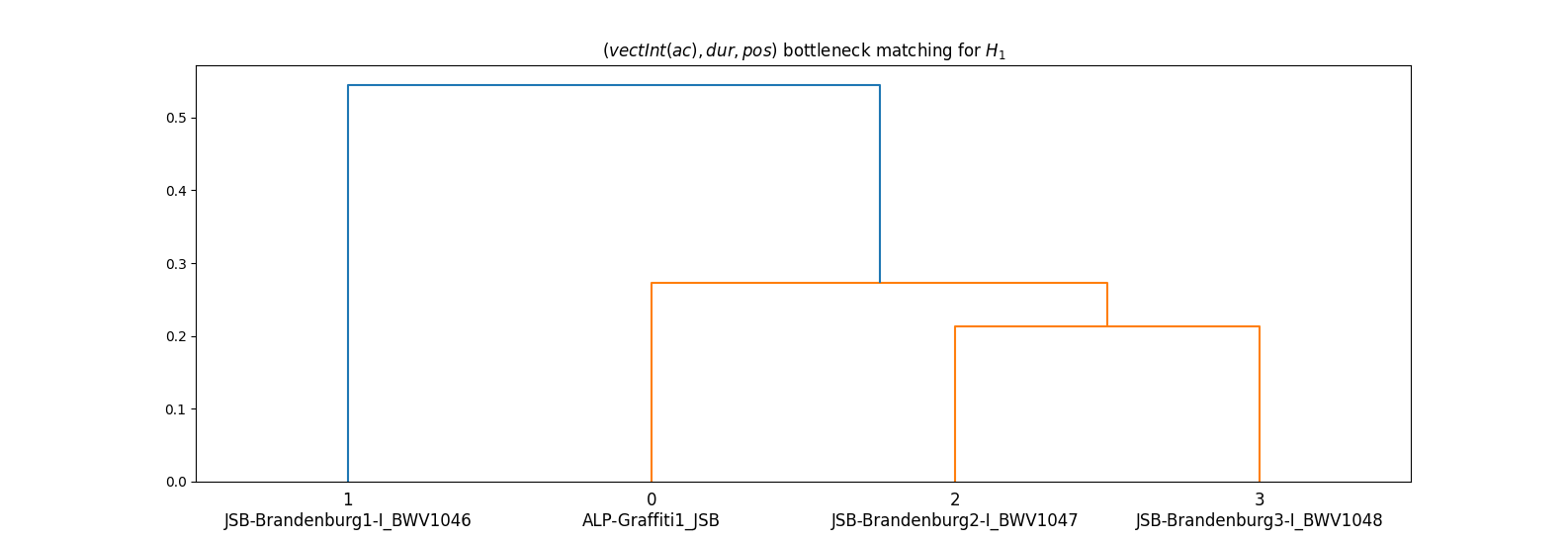}
\par\end{centering}

\protect\caption{{\scriptsize{}\label{fig:Dendrogram-for-bottleneck_H1 cat 2}Dendrogram
for bottleneck distances between $H_{1}\text{\textendash}$diagrams
for }\textbf{\scriptsize{}data mapping II}{\scriptsize{} from mm.
1\textendash 4 of Luna's }\emph{\scriptsize{}Graffiti Hommage to J.S.B.}{\scriptsize{}
and J.S. Bach's }\emph{\scriptsize{}Brandenburg Concertos BWV 1046-1048}.}
\end{figure}

\begin{figure}[btp]
\begin{centering}
\includegraphics[scale=0.37]{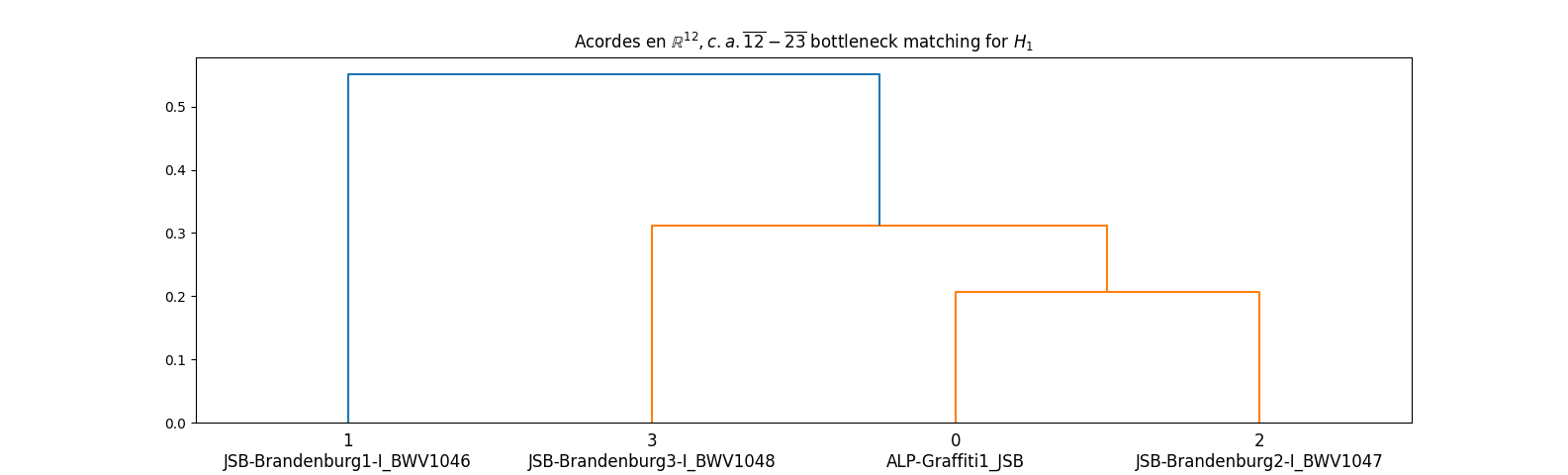}
\par\end{centering}

\protect\caption{{\scriptsize{}\label{fig:Dendrogram-for-bottleneck_H1 cat 5}Dendrogram
for bottleneck distances between $H_{1}\text{\textendash}$diagrams
for }\textbf{\scriptsize{}data mapping V}{\scriptsize{} from mm. 1\textendash 4
of Luna's }\emph{\scriptsize{}Graffiti Hommage to J.S.B.}{\scriptsize{}
and J.S. Bach's }\emph{\scriptsize{}Brandenburg Concertos BWV 1046-1048}.}
\end{figure}

\begin{figure}[btp]
\begin{centering}
\includegraphics[scale=0.37]{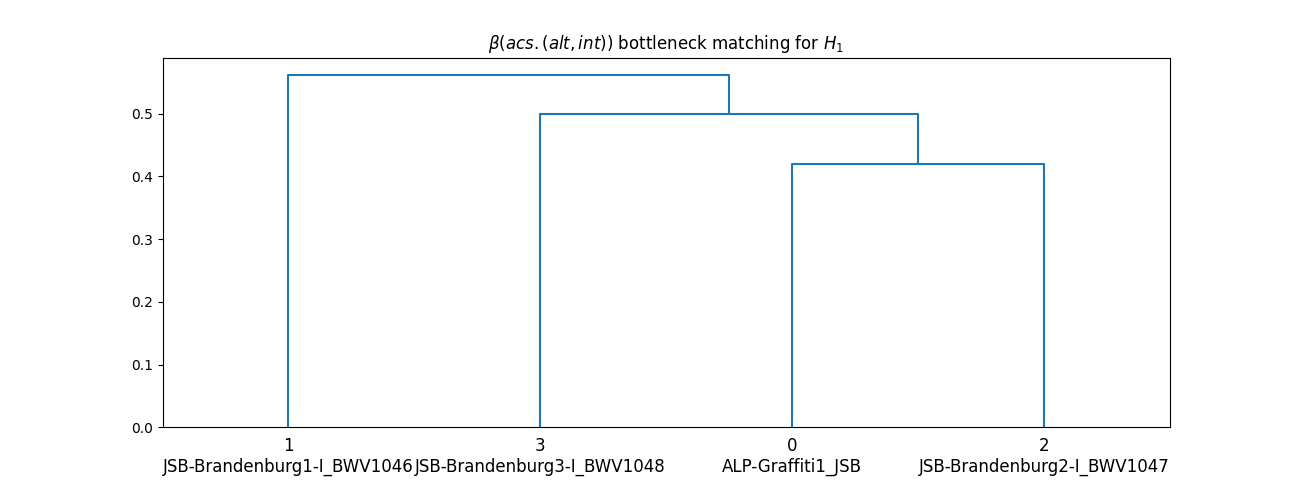}
\par\end{centering}

\protect\caption{{\scriptsize{}\label{fig:Dendrogram-for-bottleneck_H1 cat 6}Dendrogram
for bottleneck distances between $H_{1}\text{\textendash}$diagrams
for }\textbf{\scriptsize{}data mapping VI}{\scriptsize{} from mm.
1\textendash 4 of Luna's }\emph{\scriptsize{}Graffiti Hommage to J.S.B.}{\scriptsize{}
and J.S. Bach's }\emph{\scriptsize{}Brandenburg Concertos BWV 1046-1048}.}
\end{figure}

In table \ref{tab: Interpretacion dendrograma ALP-G.H. a J.S.B. vs 3 Brandenburg}
we summarize the results of our hierarchical clustering. From this
table we can conclude that \emph{Brandenburg Concerto no. 2} lies
``between'' \emph{Brandenburg Concertos nos. 1 }and\emph{ 3}, which
sometimes appear almost as opposites. Also, we observe that Luna's
\emph{Graffiti} is clearly the most distinguishable among all four
sample pieces, though it shares some common harmonic features with
the rest, as witnessed by distances between the persistence diagrams
of certain mappings: it is clustered together with \emph{Brandenburg
3} in mappings I and IV, \emph{Brandenburg 1 }in mapping III, and
\emph{Brandenburg }2\emph{ }in mapping IV (for $H_{1}\text{\textendash}$diagrams)\emph{.}

\begin{center}
\begin{table}[!tph]
\begin{centering}
{\tiny{}}%
\resizebox{\textwidth}{!}{\begin{tabular}[t]{|c|>{\centering}m{1.85cm}|>{\centering}m{1.85cm}|>{\centering}m{1.85cm}|>{\centering}m{1.85cm}|>{\centering}m{1.85cm}|>{\centering}m{1.85cm}|}
\cline{2-7} 
\multicolumn{1}{c|}{} & \multicolumn{6}{c|}{{\tiny{}\cellcolor{lightgray}}\textbf{\tiny{}Data mapping}}\tabularnewline
\cline{2-7} 
\multicolumn{1}{c|}{} & \textit{\tiny{}\cellcolor{purplepizzazz}I} & \textit{\tiny{}\cellcolor{cyan}II} & \textit{\tiny{}\cellcolor{yellow}III} & \textit{\tiny{}\cellcolor{orange}IV} & \textit{\tiny{}\cellcolor{green}V} & \textit{\tiny{}\cellcolor{blue-magenta}VI}\tabularnewline
\hline 
\multicolumn{7}{|c|}{\textit{\tiny{}\cellcolor{lightgray}}{\tiny{}$\boldsymbol{H_{0}\text{\textendash}}$}\textbf{\tiny{}diagrams}}\tabularnewline
\hline 
{\tiny{}\cellcolor{lightgray}}$<$ & \textit{\tiny{}\cellcolor{purplepizzazz}B.C.1}\\
\textit{\tiny{} B.C.2} & \textit{\tiny{}\cellcolor{cyan}B.C.2}\\
\textit{\tiny{}B.C.3} & \textit{\tiny{}\cellcolor{yellow}B.C.2}\\
\textit{\tiny{}B.C.3} & \textit{\tiny{}\cellcolor{orange}Diagrams}\\
\textit{\tiny{}from } & \textit{\tiny{}\cellcolor{green}B.C.1}\\
\textit{\tiny{}B.C.2} & \textit{\tiny{}\cellcolor{blue-magenta}B.C.1}\\
\textit{\tiny{}B.C.2}\tabularnewline
\cline{1-4} \cline{6-7} 
\multicolumn{4}{|c|}{} & \textit{\tiny{}\cellcolor{orange}all samples } & \multicolumn{2}{c|}{}\tabularnewline
\cline{1-4} \cline{6-7} 
{\tiny{}\cellcolor{lightgray}}$>$ & \textit{\tiny{}\cellcolor{purplepizzazz}B.C.3} & \textit{\tiny{}\cellcolor{cyan}G.H. to J.S.B.} & \textit{\tiny{}\cellcolor{yellow}B.C.1}\\
\textit{\tiny{}G.H. to J.S.B.} & \textit{\tiny{}\cellcolor{orange}are}\\
\textit{\tiny{}equidistant.} & \textit{\tiny{}\cellcolor{green}G.H. to J.S.B.} & \textit{\tiny{}\cellcolor{blue-magenta}G.H. to J.S.B.}\tabularnewline
\hline 
\multicolumn{7}{|c|}{\textit{\tiny{}\cellcolor{lightgray}}{\tiny{}$\boldsymbol{H_{1}\text{\textendash}}$}\textbf{\tiny{}diagrams}}\tabularnewline
\hline 
{\tiny{}\cellcolor{lightgray}}$<$ & \textit{\tiny{}\cellcolor{purplepizzazz}B.C.3}\\
\textit{\tiny{}G.H. to J.S.B.} & \textit{\tiny{}\cellcolor{cyan}B.C.2}\\
\textit{\tiny{}B.C.3} &  & \textit{\tiny{}\cellcolor{orange}B.C.1}\\
\textit{\tiny{}B.C.2} & \textit{\tiny{}\cellcolor{green}B.C.2}\\
\textit{\tiny{}G.H. to J.S.B.} & \textit{\tiny{}\cellcolor{blue-magenta}B.C.2}\\
\textit{\tiny{}G.H. to J.S.B.}\tabularnewline
\hline 
\multicolumn{6}{|c}{} & \tabularnewline
\hline 
{\tiny{}\cellcolor{lightgray}}$>$ & \textit{\tiny{}\cellcolor{purplepizzazz}B.C.2} & \textit{\tiny{}\cellcolor{cyan}B.C.1} &  & \textit{\tiny{}\cellcolor{orange}B.C.3}\\
\textit{\tiny{}G.H. to J.S.B.} & \textit{\tiny{}\cellcolor{green}B.C.1} & \textit{\tiny{}\cellcolor{blue-magenta}B.C.1}\tabularnewline
\cline{1-1} \cline{3-3} \cline{5-7} 
\end{tabular}}
\par\end{centering}{\tiny \par}

\centering{}{\footnotesize{}\protect\caption{{\scriptsize{}\label{tab: Interpretacion dendrograma ALP-G.H. a J.S.B. vs 3 Brandenburg}Closest
pairs ($<$) and most distant samples
($>$) according to the bottleneck
distance between their $H_{0}\text{\textendash}$ and $H_{1}\text{\textendash}$
persistent diagrams, by mapping}.}
}{\footnotesize \par}
\end{table}

\par\end{center}

\section*{Conclusions and future work.\label{sec:FutureWork}}

Through the study of topological features of ``clouds'' of events
represented in a music score, we get a plausible way of describing
and comparing musical features; in the case of this paper, the harmonic
structure of a fragment. We believe that the combination of the homological
descriptors of different sets of data points (mappings) associated
with a given music fragment, may lead to a homological fingerprint of
symbolic music scores, perhaps focusing specifically in some aspect,
such as harmony, or encompassing a wide range of music parameters
(rhythm, timbre, dynamics, etc.) simultaneously. Even beyond a concrete musical or musicological interpretation, by means of the (\emph{persistent}) homological invariants of musical data, evidence in other works suggests we may be able to classify and relate musical styles and features, particularly by using them as training data for machine learning models (see, for example \citet{TDAImprovesDeepLearning,UnsupervisedLearningTDA}).

Our immediate goal is to run a variety of data analysis algorithms
on results obtained from a large collection of samples, including some standard statistical analyses and recent TDA methods such as persistence landscapes (\citet{BubenikStatPersLandscapes,BubenikPersistenceToolbox,EulerCharSurfaces}), and the Euler characteristic curve (\citet{EulerCharSurfaces}). This will
give us robust mathematical descriptors for style classification, all suitable for machine learning. In the way, we will record our results under the several proposed mappings, in order to unveil which of them seem to detect certain harmonic features, or if some are disposable or redundant. We will also consider microtonal scores, to test if the present tools may lead to characterize important aspects of Indian classical music and other musical traditions.

Parting from the framework discussed in this paper, we will devote some work to the description of the dynamical properties of sequences of simplicial complexes dealt with herein.
Some work already done in studying time-varying simplicial complexes can be seen in \citet{TowTopFing,DynTopTools,HomPersTimeSeries}.
One of the ways we will approach this task is modeling harmonic progressions through dynamical systems defined on the families of simplicial complexes constructed in sections \ref{sub:Simplicial-complexes of chords by pitch}
and \ref{sub:Simplicial-complexes-of chord by pitch and interval}. Such dynamical systems may be used to study change in musical structure, and thus may lead to, for example, a general way of describing musical
change in texture and form.  On the other hand, we will incorporate a dynamical perspective on persistence diagrams and landscapes.    

Last, given that simplicial complexes introduced in this article represent the agglomeration
of chords in a given event interval, but not the way they are linked
from one to the next, certain musical aspects such as voice leading is
left out of the present analysis algorithm.  To address this issue, in a subsequent paper we will introduce two other families of simplicial complexes to model chord connections.

\pagebreak

\textcolor{red}
{
\section*{Appendix \label{sec:Appendix}}
}
For the sake of fluency, as well as to facilitate citation and reference for the reader, we devote this section to present a brief summary of simplicial and persistent homology. The basic objects of study in simplicial homology are \textbf{simplices}, which are abstract (either geometrical or merely combinatorial) equivalents of $n$-dimensional triangles (that is, vertices, edges, triangles, tetrahedrons, etc., respectively in dimensions $0,1,2,3,$ etc.).
\begin{figure}[bp]
\begin{centering}
\includegraphics[scale=0.5]{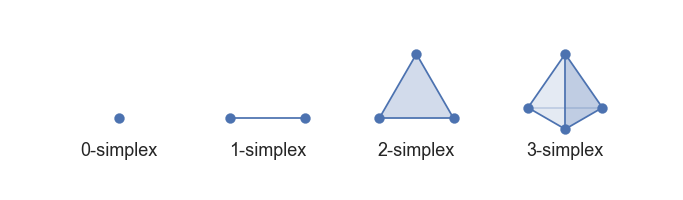}
\par\end{centering}
\protect\caption{{\label{fig: simplices examples}Geometric representation of simplices in dimensions $0,1,2$ and $3$.
}}
\end{figure}
 
   Simplices are the building blocks for the actual objects of study in this area: \textbf{simiplicial complexes}, which are topological spaces obtained by assembling simplices. \\
\begin{figure}[bp]
\begin{centering}
\includegraphics[scale=0.1]{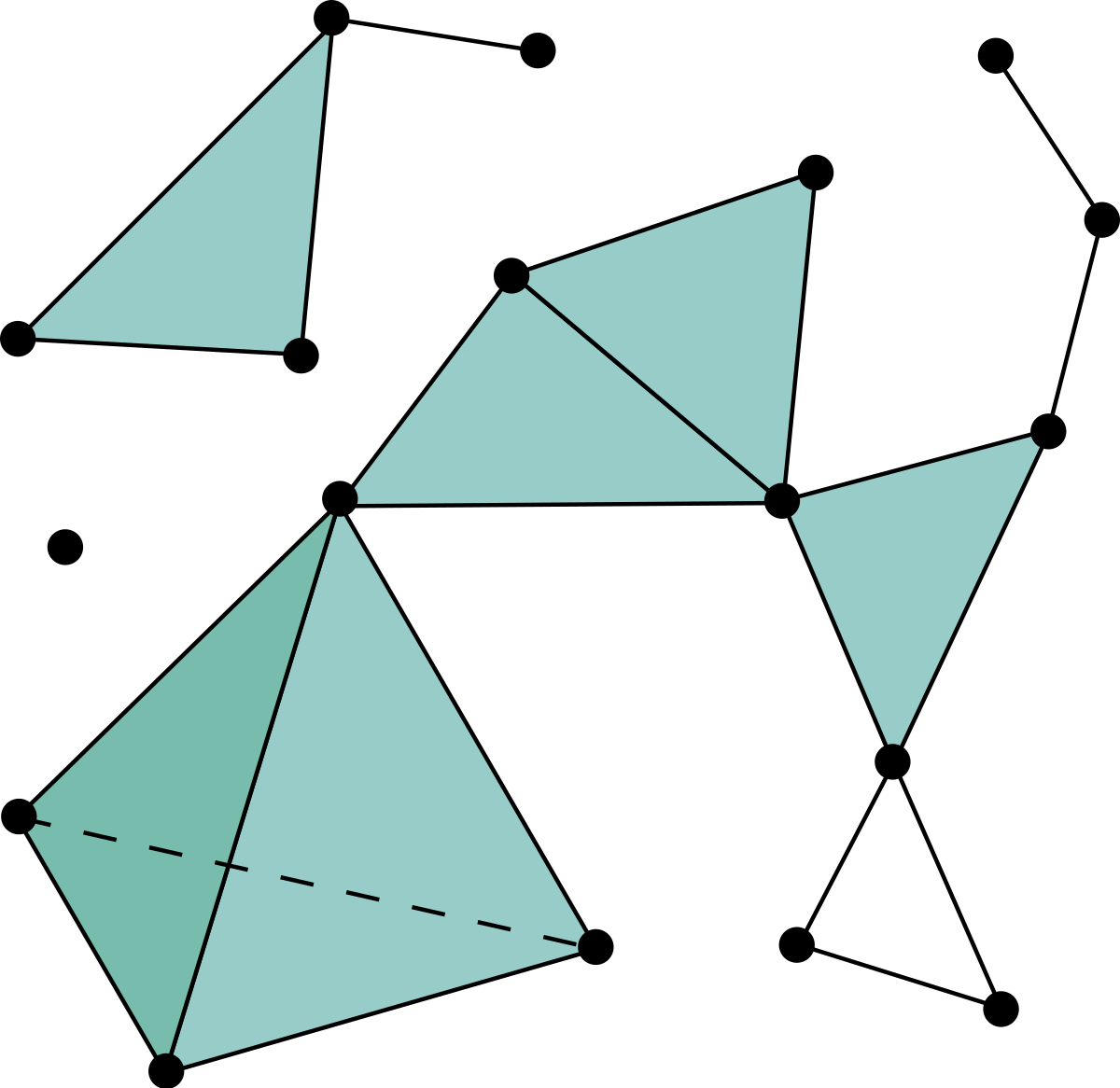}
\par\end{centering}

\protect\caption{{\label{fig: simplicial complex example} Example of a simplicial complex. }}
\end{figure}

As well as simplices, simplicial complexes can also be actual geometric objects in some Euclidean space, or more abstract, purely combinatorial objects. These two conceptions (geometric \emph{vs.} combinatorial or abstract) are equivalent, both connected by the notions of \textbf{geometric realization} and \textbf{scheme}: for every abstract $n$-dimensional simplicial complex we may build a geometric simplicial complex (its geometric realization) in $\mathbb R^{2n+1}$, and from every geometric simplicial complex we get the abstract simplicial complex (its scheme) given by the sets of vertices of its simplices. For the sake of briefness, we will be focusing on abstract simplicial complexes only. We recall that throughout the whole text, $\mathbb{R}^n$ denotes the $n$-dimensional Euclidean space, in which we consider the usual (Euclidean) distance between two points, given by the length of the straight line segment joining them. \\

Simplicial complexes are important in our context because of two main reasons: 
\begin{itemize}
\item We can encode chords (vertical events) as simplicial complexes in different ways (see section \ref{sub:Two-harmonic-simplicial-complexes}) that let us describe a fragment of a music score in terms of their homological invariants (see below). 
\item Simplicial complexes are the basis for topological data analysis, in which simplicial complexes are built (in different ways) from points representing data in a metric space, in relation to their distance (in our case, the Euclidean distance in some $\mathbb R^n $). This construction yields a filtration of simplicial complexes (in our case, the Vietoris-Rips filtration) whose vertices are these \emph{data points}. This filtration lets us have a homological description of the topological \emph{``shape''} of the given representation of the data through different \emph{``levels''}: its \textbf{persistent homology}.
\end{itemize}

Formally, we have the following (most is taken from \citet{PersHomSurvey}):

\begin{definition}
An \textbf{(abstract) simplicial complex} is a finite collection of sets $\mathcal K$ such that $\sigma \in \mathcal K$ and $\tau \subseteq \sigma$ implies $\tau \in \mathcal K$. The sets in $\mathcal K$ are referred to as its \textbf{simplices}, the union of which makes up its set of \textbf{vertices}. To explicitly refer to the vertices of a simplex we write $\sigma [x_0 , \, ... \, , x_n]$, and we say that $x_0 , \, ... \, , x_n$ \textbf{span} $\sigma$. The \textbf{dimension of a simplex} $\sigma$ is given by $|\sigma|-1$, where $|\sigma|$ denotes the cardinality of $\sigma$. A simplex on $n+1$ vertices has dimension $n$, and is called an $n$-simplex. The \textbf{dimension of a simplicial complex} is the maximum dimension of any of its simplices. A non-empty subset $\tau \subseteq \sigma$ is called a \textbf{face} of $\sigma$. Note that from the definition, a simplicial complex contains all the faces of its simplices. In the case of geometric simplicial complexes, in addition to this, simplices must be assembled together along their faces.
\end{definition}
For simplicial complexes, we have a way of algebraically encoding its shape, by describing its boundary as a formal sum of its simplices and their faces. A simplex together with a fixed order of its vertices is called an \textbf{oriented simplex}. We introduce the \textbf{boundary operator}, denoted by $\partial $, which is defined on every oriented simplex and then linearly extended to the graded abelian group (or vector space) of formal sums of simplices 
\[
\underset{ k=0 }{\overset{\text{dim} ( \mathcal{K} )} \oplus}C_n \, ,
\]
where $C_n$ is the group with basis $\mathcal S_n$, the set of simplices of dimension $n$. The boundary operator is defined as follows, for each $n$-simplex $\sigma$:

\begin{eqnarray*}
\sigma [x_0 , \, ... \, , x_n] & \overset{\partial}{\mapsto} & \underset{ k=0 }{\overset{n}{\Sigma}}(-1)^{k}\sigma[x_0 , \, ... \, , \hat x_k ,\, ... \, , x_n] \, \text{,}
\end{eqnarray*}
where $\sigma[x_0 , \, ... \, , \hat x_k ,\, ... \, , x_n]$ denotes the $(n-1)$-simplex on vertices $\{x_0 , \, ... \, ,x_n \} - \{x_k\}$. By linear extension, this defines a sequence of homomorphisms (respectively, linear functions)
\[
\xymatrix{
 \, ... \, \ar[r]^{\partial _{n+1}} & C_n \ar[r]^{\partial _{n}} & C_{n-1} \ar[r]^{\partial _{n-1}} & \, ... \,  \, \text{,}
} 
\]
which has the property that $\partial _n \circ \partial _{n+1} =0$ for every $n$. Thus, we may define the $n$-th \textbf{homology group} of simplicial complex $\mathcal K$ as the quotient group (resp. vector space)
\[
H_n (\mathcal K)=Ker (\partial _n) / Im (\partial _{n+1}) \, \text{.}
\]

Homology groups somehow capture the way simplices are ``glued'' together to form a simplicial complex. They give a rough description of the shape of a space by measuring how many $n$-dimensional ``holes'' or ``voids'' are enclosed by it. This is summarized by the ranks or dimensions of the homology groups, called the  \textbf{Betti numbers} of $\mathcal K$, and denoted by $\beta _0 (\mathcal K) , \beta _1 (\mathcal K) , \beta _2 (\mathcal K) , \, ... \,$ These numbers constitute a family of important homotopical invariants of a topological space, and are summarized in the \textbf{Euler characteristic}, given by their alternated sum:

\[
E(\mathcal K ) = \underset{ k=0 }{\overset{\text{dim}(\mathcal K)}{\Sigma}}(-1)^{k}\beta _k (\mathcal K) \, \text{.}
\]

Betti numbers are central in topological data analysis, particularly in persistent homology. \textbf{Persistent homology} is the computation of homological features (namely, Betti numbers) at different ``levels'', ``scales'' or ``resolutions'' to get an algebraic description of the shape of a set of points $S$ in a space. This leads to consider filtrations (another crucial concept in persistent homology). A \textbf{filtration} of a simplicial complex $\mathcal K$ is a sequence of simplicial subcomplexes 
\[
F_0 \subset F_1 \subset F_2 \subset \, ... \, \subset F_{\ell}=\mathcal K \, \text{.}
\]
We restrict ourselves to finite filtrations.

We consider a one-paremeter filtration of simplicial complexes built upon the points in $S$ (taken as vertices). This parameter, usually $\varepsilon$, establishes the level or scale of scope on the data. 

\[
S=F_{\varepsilon_0} \subset F_{\varepsilon_1} \subset F_{\varepsilon_2} \subset \, ... \, \subset F_{\varepsilon_{\ell}}=\mathcal K \, \text{,}
\]
where $\mathcal K$ is the simplex on all points in $S$ and $\varepsilon_0 <\varepsilon_1  <\, ... \, <\varepsilon_{\ell}$. Simplicial complexes in the filtration can be defined according to different constructions. There are different ways to associate a filtration of simplicial complexes with a set of discrete points in a metric space. Persistent homology consists of computing the Betti numbers (which \emph{count} $n$-dimensional holes or voids) of the simplicial complexes in the filtration considered.

Now we present one particular filtration that can be associated with any cloud of data points in a metric space, the Vietoris-Rips filtration, one of the most commonly used for computing persistent homology, and the one incorporated by the Python library used in this work.

\begin{figure}[bp]
\begin{centering}
\includegraphics[scale=0.5]{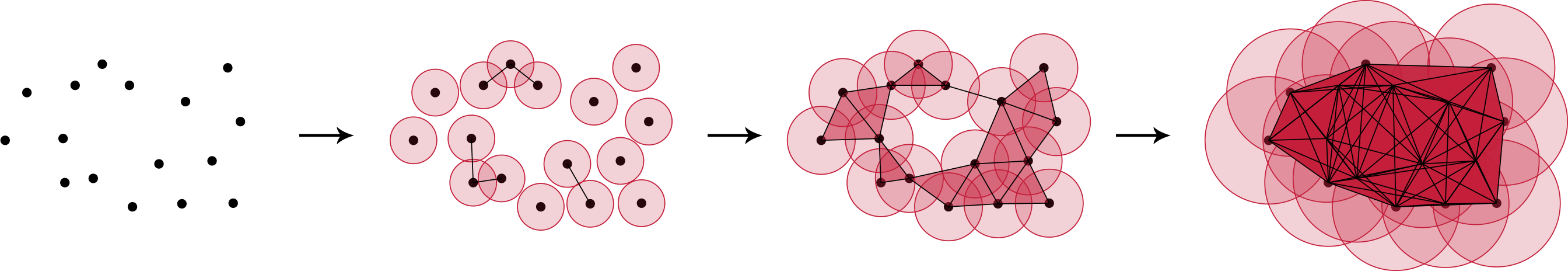}
\par\end{centering}

\protect\caption{{\label{fig: Vietoris-Rips example} Vietoris-Rips filtration of a point set at four different values of $\varepsilon$. }}
\end{figure}

Given a metric space $(X,d)$ and a finite set of points $S=\{x_0 , x_1 , \, ... \, , x_N\}$ in X, for each given distance $\varepsilon$ we consider the Vietoris-Rips simplicial complex 
\[
\mathcal {VR} (S,\varepsilon) = \{ \sigma \subseteq S \mid \text{diam} (\sigma) \leq 2\varepsilon \} \, \text{.}
\]
$\mathcal {VR} (S,\varepsilon)$ is the simplicial complex whose simplices consist  of those subsets of points in $S$ which are not further than $2\varepsilon$ among themselves. So simplices in this case represent ``closeness'' of points: a $k$-simplex is formed whenever $k+1$ points can be enclosed together in a $k$-ball of radius $\varepsilon$. As the parameter $\varepsilon$ varies, we obtain a filtration of simplicial complexes, called the \textbf{Vietoris-Rips filtration}, in which the last element is the $N$-dimensional simplex on all vertices of $S$. Figure \ref{fig: Vietoris-Rips example} illustrates this definition.

The Betti numbers of the simplicial complexes in this filtration are computed to obtain a homological description of the ``shape'' of our cloud of points $S$ in space $(X,d)$ at different scales.  They are usually summarized in two graphic ways: persistent barcode graphs and persistent diagrams. These barcodes and diagrams (see figures in section \ref{sub: PersHomMaps}) record the changing values of the Betti numbers throughout the Vietoris-Rips filtration, and so give us a way of visualizing at which scales do connected components and $n$-dimensional voids appear (their birth, $b$) and disappear (their death, $d$). Barcodes are plotted as a set of line segments in the plane $\epsilon \times \text{Betti numbers}$, starting and ending at the values of $\epsilon$ for which each homological feature persists. Similarly, persistence diagrams show the points $(b,d)$ in the plane $\mathbb{R}^2 = Birth \times Death$, usually plotted alongside the diagonal (all points sit above this line). From this information, certain conclusions about the general distribution of the given points can be drawn, which help in understanding and classifying large collections of data. For example, the longest barcodes are interpreted as the most relevant (persistent) features of the point cloud's shape.

\begin{figure}

     \begin{minipage}[t]{0.45\textwidth}

         \includegraphics[width=1\textwidth]{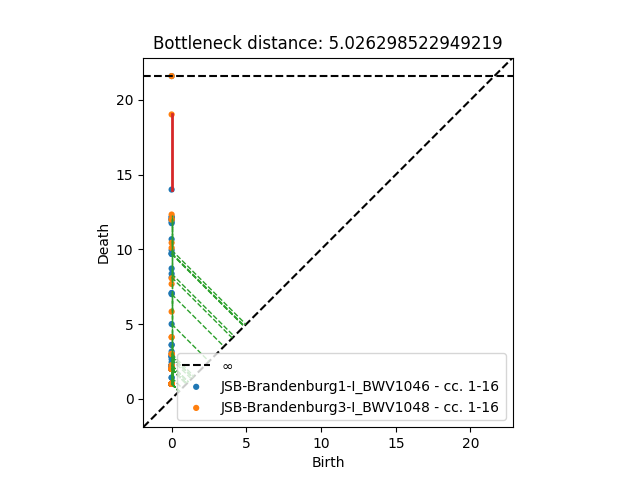}
         %\caption{\tiny Matching for the bottleneck distance between persistence diagrams in dimension $0$ of J. S. Bach's Brandenburg Concertos 1 and 3.}
         \label{fig: bottleneck_dist_JSB-Brandenburg1-I-cc1-16-vs-JSB-Brandenburg3-II-cc1-16-TDA-III-H0}
     \end{minipage}
     \begin{minipage}[t]{0.45\textwidth}

         \includegraphics[width=1\textwidth]{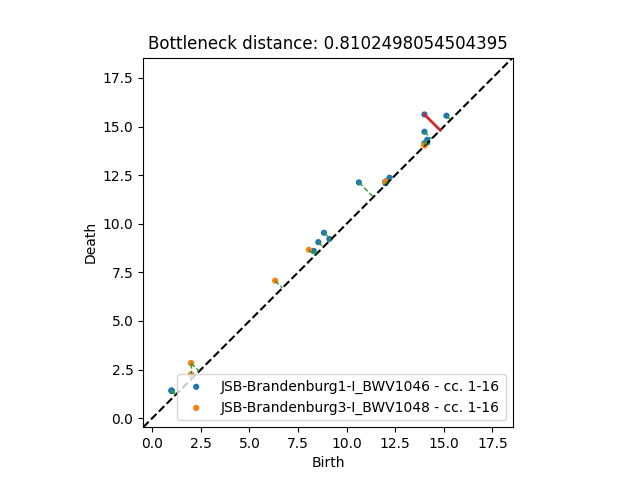}
         %\caption{\tiny Matching for the bottleneck distance between persistence diagrams in dimension $1$ of J. S. Bach's Brandenburg Concertos 1 and 3.}
         \label{bottleneck_dist_JSB-Brandenburg1-I-cc1-16-vs-JSB-Brandenburg3-II-cc1-16-TDA-III-H1}
     \end{minipage}
     \caption{\label{fig: Matching for the bottleneck distance JSB Brandenburg 1-3 mapping V} \tiny Matching for the bottleneck distance between persistence diagrams in dimensions $0$ (left) and $1$ (right) for measures 1-16 of J. S. Bach's Brandenburg Concertos 1 and 3, under mapping V.}
\end{figure}

Barcodes and persistence diagrams corresponding to different sets of points can be compared in several ways. One of the most common is the \textbf{bottleneck distance}, usually denoted $W_{\infty}$ (in relation to the family of Wasserstein metrics). This consists of finding a  pairing between points $x=(b,d)$ in both diagrams, minimizing the maximum possible  $L_\infty$-distance between corresponding points in the pairing. Formally, the bottleneck distance between diagrams $D,D'$ is given by 
\[
W_{\infty} (D,D') = \underset{\varphi : D \leftrightarrow D'}{\text{inf}} \underset{x \in D}{\text{sup}} \, \, \lVert x-\varphi (x) \rVert _{\infty} \, \text{,}
\]
where $\varphi$ ranges over all bijections between $D$ and $D'$, and
\[
\lVert (x_1 , x_2)-(y_1 , y_2) \rVert _{\infty} = \text{max} \{ \mid x_1 - y_1\mid , \mid x_2 - y_2 \mid \} \, \text{.}
\]
Whenever there is no bijection between $D$ and $D'$, a partial pairing is considered, and those remaining points are paired with their projection over the diagonal. This metric has the property that it is the least number such that one can draw on the plane squares of side $2W_{\infty} (D,D') $ centered at the elements of $D$, and these will also contain the corresponding points of $D'$ under the matching that defines $W_{\infty} (D,D') $. Figure \ref{fig: Matching for the bottleneck distance JSB Brandenburg 1-3 mapping V} shows how points in two different persistence diagrams are matched in order to compute their bottleneck distance. Distances taken among several persistence diagrams can be plotted in a dendrogram  (see section \ref{sec:Results.}), which facilitates visualization and comparison of the values obtained for several samples.

\begin{figure}

     \begin{minipage}[t]{0.45\textwidth}

         \includegraphics[width=.85\textwidth]{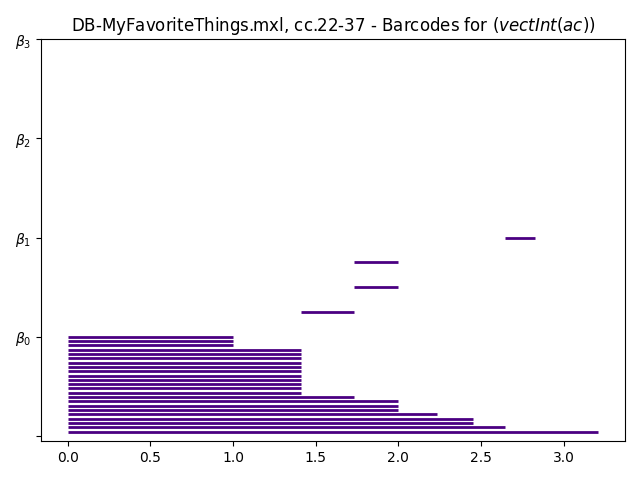}
         \caption*{\tiny Persistent homology barcodes for interval vectors.}
         \label{fig: DB-MyFavoriteThings-cc22-37-TDA-I-barcodes}
     \end{minipage}
     \begin{minipage}[t]{0.45\textwidth}

         \includegraphics[width=.85\textwidth]{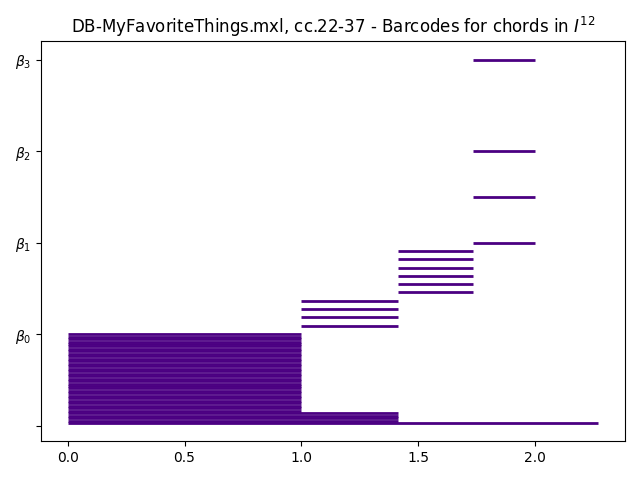}
         \caption*{\tiny Persistent homology barcodes for binary pitch vectors.}
         \label{DB-MyFavoriteThings-cc22-37-TDA-II-barcodes}
     \end{minipage}
\medskip

     \begin{minipage}[t]{0.45\textwidth}

         \includegraphics[width=.85\textwidth]{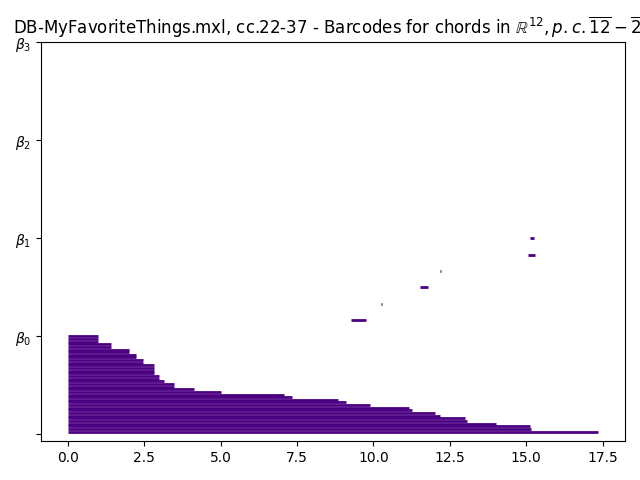}
         \caption*{\tiny Persistent homology barcodes for embedded normal form vectors.}
         \label{fig: DB-MyFavoriteThings-cc22-37-TDA-III-barcodes}
     \end{minipage}
     \begin{minipage}[t]{0.45\textwidth}

         \includegraphics[width=.85\textwidth]{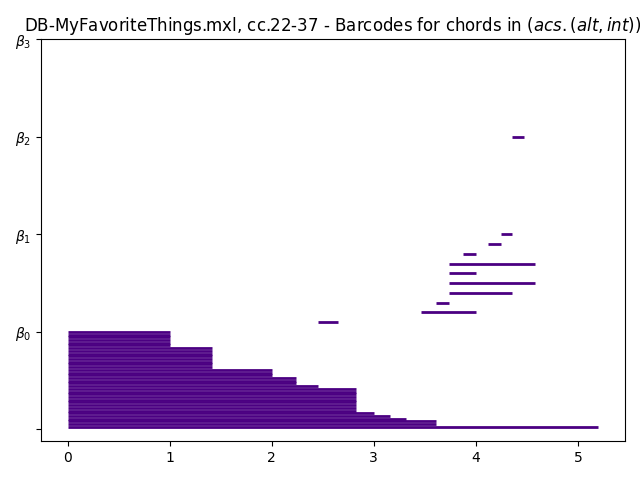}
         \caption*{\tiny Persistent homology barcodes for pitch-interval vectors.}
         \label{fig: DB-MyFavoriteThings-cc22-37-TDA-IV-barcodes}
     \end{minipage}
\caption{\scriptsize Barcodes for the first phrase in Dave Brubeck's version of 'My favourite things', by Rodgers and Hammerstein.}
        \label{fig: TDA - DB-MyFavoriteThings-cc22-37}
\end{figure}

\begin{figure}

     \begin{minipage}[tp]{0.45\textwidth}

         \includegraphics[width=.85\textwidth]{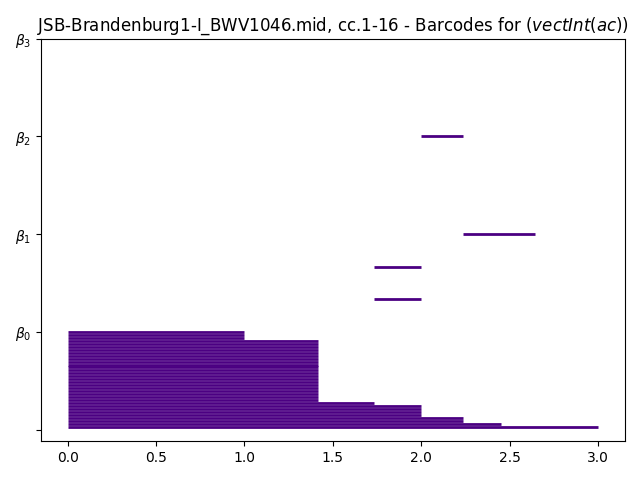}
         \caption*{\tiny Persistent homology barcodes for interval vectors.}
         \label{fig: DB-MyFavoriteThings-cc22-37-TDA-I-barcodes}
     \end{minipage}
     \begin{minipage}[tp]{0.45\textwidth}

         \includegraphics[width=.85\textwidth]{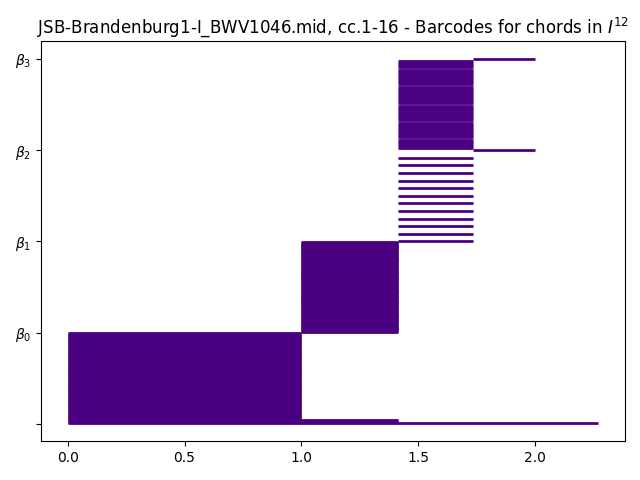}
         \caption*{\tiny Persistent homology barcodes for binary pitch vectors.}
         \label{DB-MyFavoriteThings-cc22-37-TDA-II-barcodes}
     \end{minipage}
\medskip

     \begin{minipage}[tp]{0.45\textwidth}

         \includegraphics[width=.85\textwidth]{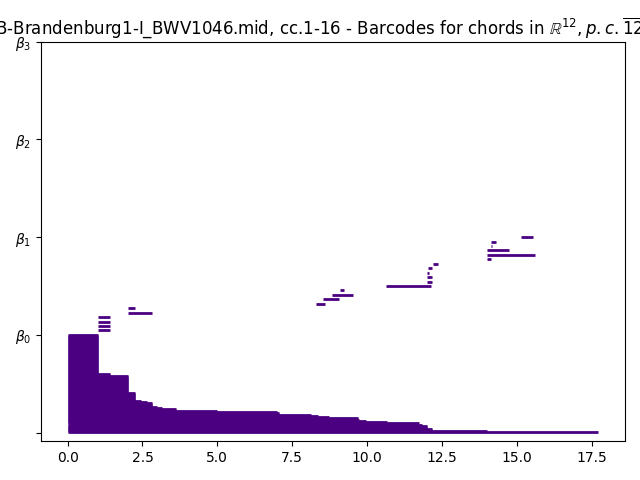}
         \caption*{\tiny Persistent homology barcodes for embedded normal form vectors.}
         \label{fig: DB-MyFavoriteThings-cc22-37-TDA-III-barcodes}
     \end{minipage}
     \begin{minipage}[tp]{0.45\textwidth}

         \includegraphics[width=.85\textwidth]{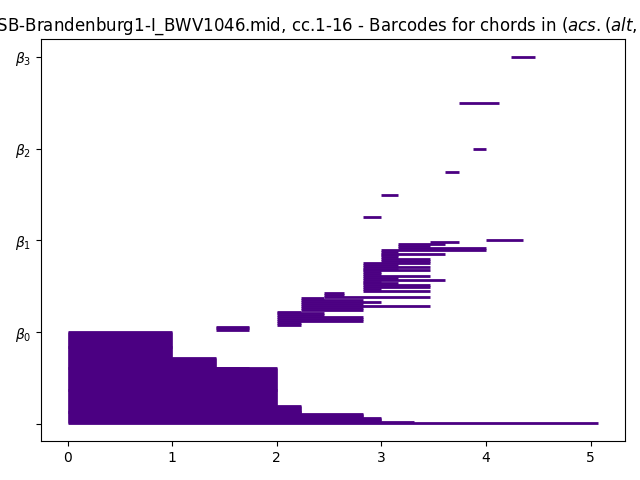}
         \caption*{\tiny Persistent homology barcodes for pitch-interval vectors.}
         \label{fig: DB-MyFavoriteThings-cc22-37-TDA-IV-barcodes}
     \end{minipage}
        \caption{\scriptsize Barcodes for mm. 1-16 from J. S. Bach's Brandenburg Concerto no. 1 - I.}
        \label{fig: TDA - JSB-Brandenburg1-I-cc1-16}
\end{figure}

\begin{figure}

     \begin{minipage}[tp]{0.45\textwidth}

         \includegraphics[width=.85\textwidth]{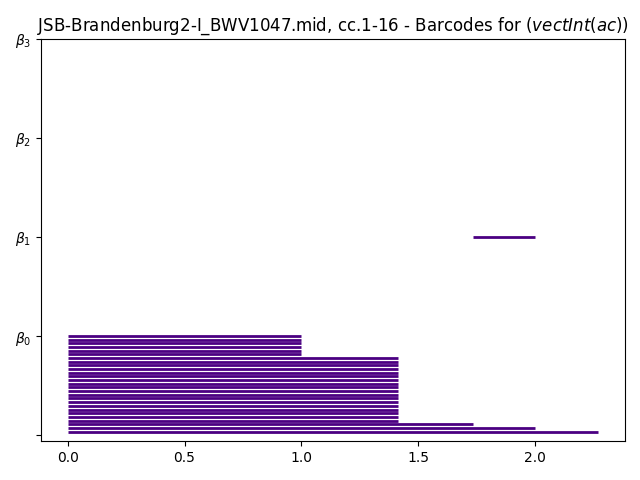}
         \caption*{\tiny Persistent homology barcodes for interval vectors.}
         \label{fig: DB-MyFavoriteThings-cc22-37-TDA-I-barcodes}
     \end{minipage}
     \begin{minipage}[tp]{0.45\textwidth}

         \includegraphics[width=.85\textwidth]{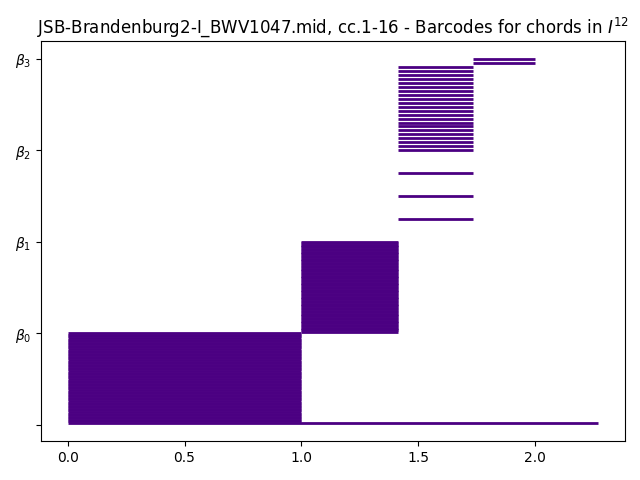}
         \caption*{\tiny Persistent homology barcodes for binary pitch vectors.}
         \label{DB-MyFavoriteThings-cc22-37-TDA-II-barcodes}
     \end{minipage}
\medskip

     \begin{minipage}[tp]{0.45\textwidth}

         \includegraphics[width=.85\textwidth]{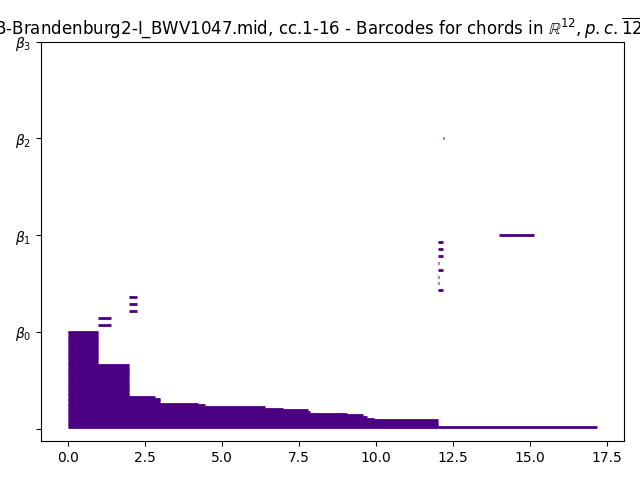}
         \caption*{\tiny Persistent homology barcodes for embedded normal form vectors.}
         \label{fig: DB-MyFavoriteThings-cc22-37-TDA-III-barcodes}
     \end{minipage}
     \begin{minipage}[tp]{0.45\textwidth}

         \includegraphics[width=.85\textwidth]{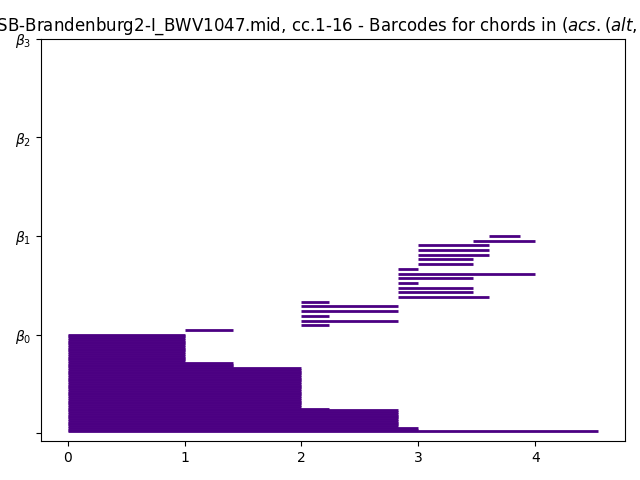}
         \caption*{\tiny Persistent homology barcodes for pitch-interval vectors.}
         \label{fig: DB-MyFavoriteThings-cc22-37-TDA-IV-barcodes}
     \end{minipage}
        \caption{\scriptsize Barcodes for mm. 1-16 from J. S. Bach's Brandenburg Concerto no. 2 - I.}
        \label{fig: TDA - JSB-Brandenburg2-I-cc1-16}
\end{figure}

\begin{figure}

     \begin{minipage}[tp]{0.45\textwidth}

         \includegraphics[width=.85\textwidth]{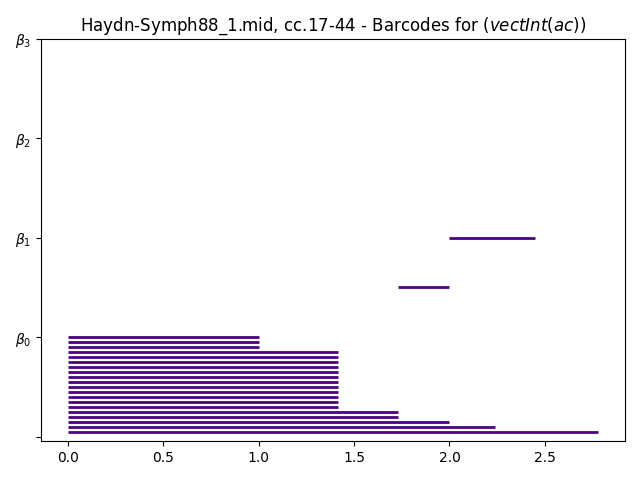}
         \caption*{\tiny Persistent homology barcodes for interval vectors.}
         \label{fig: DB-MyFavoriteThings-cc22-37-TDA-I-barcodes}
     \end{minipage}
     \begin{minipage}[tp]{0.45\textwidth}

         \includegraphics[width=.85\textwidth]{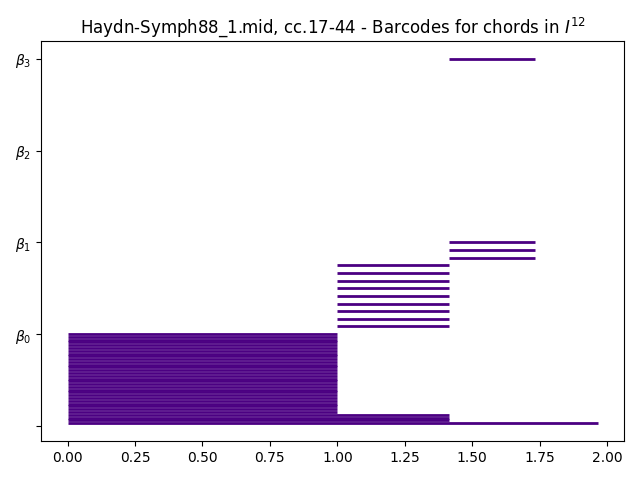}
         \caption*{\tiny Persistent homology barcodes for binary pitch vectors.}
         \label{DB-MyFavoriteThings-cc22-37-TDA-II-barcodes}
     \end{minipage}
\medskip

     \begin{minipage}[tp]{0.45\textwidth}

         \includegraphics[width=.85\textwidth]{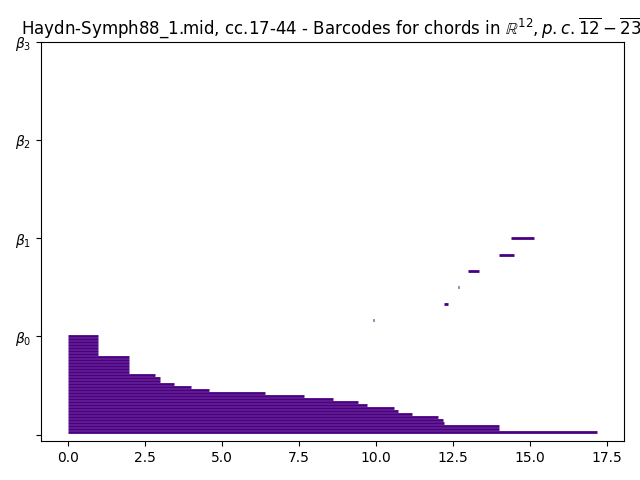}
         \caption*{\tiny Persistent homology barcodes for embedded normal form vectors.}
         \label{fig: DB-MyFavoriteThings-cc22-37-TDA-III-barcodes}
     \end{minipage}
     \begin{minipage}[tp]{0.45\textwidth}

         \includegraphics[width=.85\textwidth]{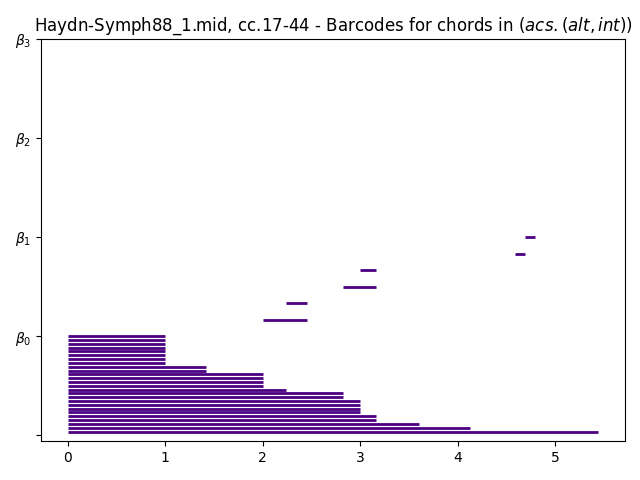}
         \caption*{\tiny Persistent homology barcodes for pitch-interval vectors.}
         \label{fig: DB-MyFavoriteThings-cc22-37-TDA-IV-barcodes}
     \end{minipage}
        \caption{\scriptsize Barcodes for mappings III-VI of mm. 17-44 from F. J. Haydn's Symphony no. 88, I. Allegro.}
        \label{fig: TDA - Haydn-Symph88_1-cc17-44}
\end{figure}

\begin{figure}

     \begin{minipage}[tp]{0.45\textwidth}

         \includegraphics[width=.85\textwidth]{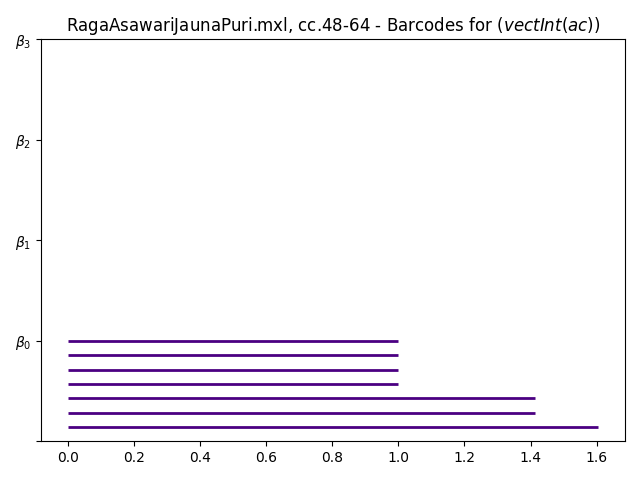}
         \caption*{\tiny Persistent homology barcodes for interval vectors.}
         \label{fig: DB-MyFavoriteThings-cc22-37-TDA-I-barcodes}
     \end{minipage}
     \begin{minipage}[tp]{0.45\textwidth}

         \includegraphics[width=.85\textwidth]{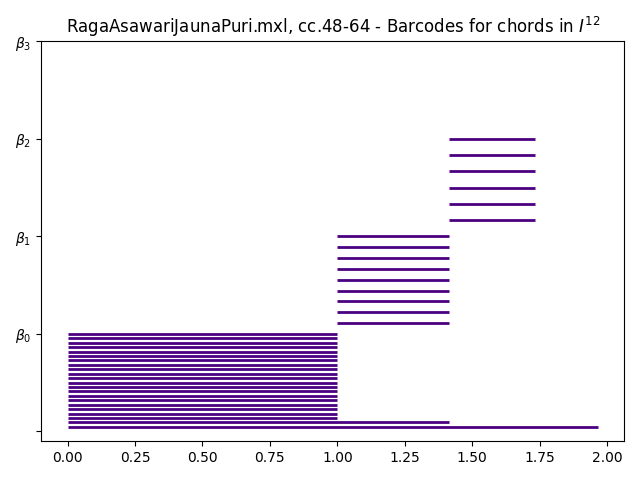}
         \caption*{\tiny Persistent homology barcodes for binary pitch vectors.}
         \label{DB-MyFavoriteThings-cc22-37-TDA-II-barcodes}
     \end{minipage}
\medskip

     \begin{minipage}[tp]{0.45\textwidth}

         \includegraphics[width=.85\textwidth]{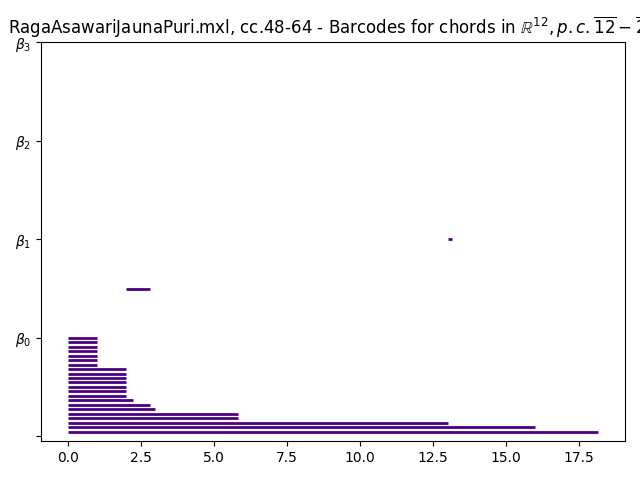}
         \caption*{\tiny Persistent homology barcodes for embedded normal form vectors.}
         \label{fig: DB-MyFavoriteThings-cc22-37-TDA-III-barcodes}
     \end{minipage}
     \begin{minipage}[tp]{0.45\textwidth}

         \includegraphics[width=.85\textwidth]{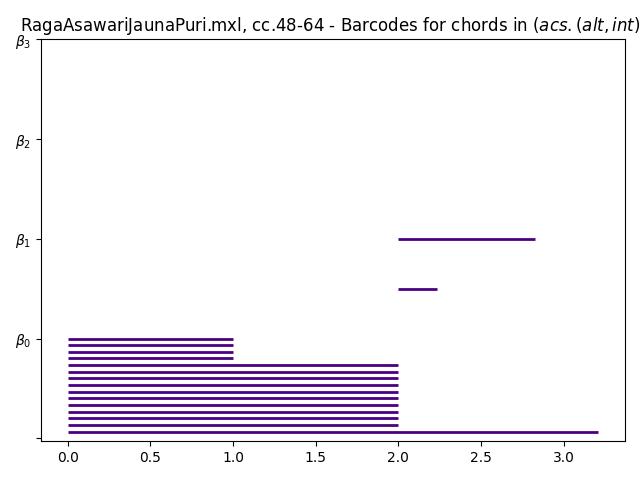}
         \caption*{\tiny Persistent homology barcodes for pitch-interval vectors.}
         \label{fig: DB-MyFavoriteThings-cc22-37-TDA-IV-barcodes}
     \end{minipage}
        \caption{\scriptsize Barcodes for mappings III-VI of a single phrase from Raga Asavari.}
        \label{fig: TDA - RagaAsawari-cc48-64}
\end{figure}

Finally, we include the plots of barcodes for mappings III-VI (only harmonic data), corresponding to fragments from several scores (figures \ref{fig: TDA - DB-MyFavoriteThings-cc22-37}-\ref{fig: TDA - RagaAsawari-cc48-64}). From looking at the illustrations we can tell there are some distinctive features in the barcode diagrams of samples from different music styles. We can also see the difference in the overall shape of barcodes resulting from mappings III and IV (forming straight blocks), against barcodes corresponding to mappings V and VI (irregular staggered and jagged shapes). In addition, we point out how barcodes for mappings III and V seem to simplify features arising from mappings IV and VI, respectively. In-depth comparison and analysis of particular examples belonging to diverse genres, styles or traditions will be discussed thoroughly in future works.

\center

\pagebreak{}

%\bibliographystyle{authordate1}
%\phantomsection\addcontentsline{toc}{section}{\refname}\bibliography{/Users/alberto/Documents/doc_biblio/docBibTeX}

No potential conflict of interest was reported by the authors. 

\section*{ORCID}% The * makes this section unnumbered.
\addcontentsline{toc}{section}{ORCID}% Do not change this.

Do not change this. Production will take care of it if the paper is accepted.

\bibliographystyle{tMAM}% The file tMAM.bst from author package must be in your working directory.
% It automatically formats the bibliography, don't change it at all. Do not include the extension .bst in this command, nor the extension .bib in the next command.
\bibliography{/Users/alberto/Documents/doc_biblio/docBibTeX}

\addcontentsline{toc}{section}{References}

\end{document}